\documentclass[twoside, 12pt]{article}
\usepackage{amssymb,mathrsfs,psfrag,graphicx}
\usepackage{hyperref}
\def\v{\varepsilon}

\def\x{\xi}
\def\t{\theta}
\def\T{\Theta}
\def\k{\kappa}

\def\m{\mu}
\def\a{\alpha}
\def\b{\beta}
\def\g{\gamma}
\def\d{\delta}
\def\l{\lambda}
\def\f{\frac}
\def\p{\phi}
\def\r{\rho}
\def\ra{\rightarrow}
\def\s{\sigma}
\def\z{\zeta}

\def\di{\displaystyle}
\def\i{\infty}

\begin{document}

\title{\bf Stability of Wave Patterns to the Inflow Problem of Full Compressible Navier-Stokes Equations} \vskip 0.5cm
\author{\bf  Xiaohong Qin\thanks{X. Qin is supported in part by  NSFC-NSAF grant
(No. 10676037). E-mail: xqin@amss.ac.cn.} \qquad  Yi Wang\thanks{Y.
Wang is supported by NSFC grant (No. 10801128) and the Knowledge
Innovation Program of the Chinese Academy of Sciences. E-mail:
wangyi@amss.ac.cn.}}
\date{} \maketitle
\noindent{\small  $^\ast$Department of Mathematics, Nanjing
University of Science and Technology, Nanjing, China

\noindent $^{\dag}$Institute of Applied Mathematics, Academy of
Mathematics and Systems Science, CAS, Beijing 100190, China}

\date{}
\maketitle

\begin{abstract}
The inflow problem of full compressible Navier-Stokes equations is
considered on the half line $(0,+\i)$. Firstly, we give the
existence (or non-existence) of the boundary layer solution to the
inflow problem when the right end state $(\r_+,u_+,\t_+)$ belongs to
the subsonic, transonic and supersonic regions respectively. Then
the asymptotic stability of not only the single contact wave but
also the superposition of the boundary layer solution, the contact
wave and the rarefaction wave to the inflow problem are investigated
under some smallness conditions. Note that the amplitude of the
rarefaction wave can be arbitrarily large. The proofs are given by
the elementary energy method.
\end{abstract}

\renewcommand{\theequation}{\arabic{section}.\arabic{equation}}
\setcounter{equation}{0}
\section{Introduction}
In this paper, we consider the half space problem of the full (or
non-isentropic) compressible Navie-Stokes equations in Eulerian
coordinate:
\begin{equation}
\left\{
\begin{array}{ll}
\di \r_t+(\r u)_x=0,&x>0,~~t>0,\\[1mm]
\di (\r u)_t+(\r u^2+p)_x=\m u_{xx},& x>0,~~t>0,\\[1mm]
\di {[\r(e+\f{u^2}{2})]}_t+[\r u(e+\f{u^2}{2})+pu]_x=\k \t_{xx}+\m
(uu_x)_x,&x>0,~~t>0,\\[1mm]
\end{array}
\right. \label{(1.1)}
\end{equation}
where $\rho(t,x)>0$ is the density, $u(t,x)$ is the velocity,
 $\theta(t,x)$ is the absolute temperature of the gas, and
 $p=p(\rho,\theta)$ is the pressure, $e=e(\rho,\theta)$ is
the internal energy, $\mu>0$ is the viscosity constant, and
$\kappa>0$ is the coefficient of heat conduction. Here we consider
the perfect gas, that is
\begin{equation}
p=R\rho\theta=A\rho^\g \exp{(\f{\g-1}{R}s)},\quad
e=\frac{R\theta}{\gamma-1}+const.,\label{(1.2)}
\end{equation}
where $s$ is the entropy, $\gamma>1$ is the adiabatic exponent, and
$A, R>0$ are gas constants.

The initial values are given by
\begin{equation}
(\r,u,\t)(t=0,x)=(\r_0,u_0,\t_0)(x)\ra(\r_+,u_+,\t_+), x\ra +\i,
\label{(1.3)}
\end{equation}
where $(\r_+,u_+,\t_+)$ is a constant state with $\r_+,\t_+$
positive. The boundary values are the following:
\begin{equation}
(\r,u,\t)(t,x=0)=(\r_-,u_-,\t_-),\label{(1.4)}
\end{equation}
where $\r_->0,~\t_->0,~u_-$ are given constants. And of course the
initial values (\ref{(1.3)}) and the boundary condition
(\ref{(1.4)}) satisfy the compatible condition at the origin
$(0,0)$.

 According to the sign of the velocity $u_-$ on the boundary $\{x=0\}$,
the following three types of problems are proposed
\cite{[Matsumura]}:

(1) the inflow problem, i.e., the velocity $u_->0$;

(2) the
outflow problem with $u_-<0$;

(3) the impermeable wall problem, i.e., $u_-=0$.

\noindent It should be remarked that in the cases (2) and (3), the
density $\r_-$ can not be given on the boundary by the  theory of
well-posedness on the hyperbolic equation $(\ref{(1.1)})_1$.

In this paper, we are interested in the case of the inflow problem
(\ref{(1.1)}), (\ref{(1.3)})-(\ref{(1.4)}). When $\k=\mu=0,$ the
compressible system (\ref{(1.1)}) becomes the inviscid Euler system
\begin{equation}
\left\{
\begin{array}{l}
\di \r_t+(\r u)_x=0,\\[1mm]
\di (\r u)_t+(\r u^2+p)_x=0,\\[1mm]
\di {[\r(e+\f{u^2}{2})]}_t+[\r u(e+\f{u^2}{2})+pu]_x=0.
\end{array}
\right. \label{(1.4+)}
\end{equation}
 The Euler system (\ref{(1.4+)}) is a typical example of the
hyperbolic conservation laws. It is well-known that the main feature
of the solutions to the hyperbolic conservation laws is the
formation of the shock wave no matter how smooth the initial values
are. The Euler system (\ref{(1.4+)}) contains three basic wave
patterns in the solutions to the Riemann problem. They are two
nonlinear waves, called shock wave and rarefaction wave, and one
linear wave called contact discontinuity. The above three dilation
invariant wave solutions and their linear superpositions in the
increasing order of characteristic speed, i.e., Riemann solutions,
govern both local and large-time behavior of solutions to the Euler
system. The invscid Euler system (\ref{(1.4+)}) is an ideal model in
gas dynamics when the dissipation effects are neglected, thus it is
of great importance to study the corresponding viscous system
(\ref{(1.1)}).

There has been a large literature on the large-time behavior of the
solutions to Cauchy problem of the compressible Navier-Stokes
equations (\ref{(1.1)}) toward the viscous versions of the three
basic wave patterns. In 1985, Matsumura-Nishihara
\cite{[Matsumura-Nishihara-1]} firstly proved the stability of the
viscous shock wave to the isentropic compressible Navier-Stokes
equations (i.e., the entropy $s$ is assumed to be constant and the
energy conservation law is not considered). Since then, many authors
had been attracted to study the stability of the viscous wave
patterns and much progress has been made. We refer to
\cite{[Huang-Li-Matsumura]}, \cite{[Huang-Matsumura-Xin-1]},
\cite{[Huang-Matsumura-Xin-2]}, \cite{[Huang-Xin-Yang]},
\cite{[Liu]}, \cite{[Liu-Xin-1]}, \cite{Liu-Xin-2},
\cite{[Kawashima-Matsumura]},
\cite{[Kawashima-Matstumura-Nishihara]},
\cite{[Matsumura-Nishihara-2]}, \cite{[Nishihara-Yang-Zhao]},
\cite{[Szepessy-Xin]}, \cite{[Xin]} and some references therein. All
these results show  that the large-time behavior  of the solutions
to Cauchy problem are basically governed by the Riemann problem of
the corresponding Euler equations.

Recently, the initial-boundary value problem (IBVP) of (\ref{(1.1)})
attracts
 increasing interest  because it has more physical meanings and of
course produces some new mathematical difficulties due to the
boundary effect. Not only basic wave patterns but also a new wave,
which is called  boundary layer solution (BL-solution for brevity)
\cite{[Matsumura]}, may appear in the IVBP case. Matsumura
\cite{[Matsumura]} proposes a criterion on the question when the
BL-solution forms to the isentropic Navier-Stokes equations. The
argument is also true to the full  Navier-Stokes equations
(\ref{(1.1)}). Consider the Riemann problem  to the Euler equations
(\ref{(1.4+)}), where the initial right  end state is given by the
far field state $(\r_+,u_+,\t_+)$ in (\ref{(1.3)}), and the left
end state $(\r_-,u_-,\t_-)$ is given by the all possible states
which are consistent with the boundary condition (\ref{(1.4)}) at
$\{x=0\}$. Note that for the outflow problem, $\r_-$ can not be
prescribed and is free on the boundary. On the one hand, when the
left end state is uniquely determined so that the value at the
boundary $\{x=0\}$ of the solution to the Riemann problem is
consistent with the boundary condition, we expect no BL-solution
occurs. On the other hand, if the value of the Riemann problem's
solution on the boundary is not consistent with the boundary
condition for any admissible left end  state, we expect a
BL-solution which compensates the gap comes up. Such BL-solution
could be constructed by the stationary solution to Navier-Stokes
equations. The existence and stability of the BL-solution (to the
inflow or outflow problems, to the isentropic or full Navier-Stokes
equations) are studied extensively by many authors, see
\cite{[Huang-Li-Shi]}, \cite{[Huang-Matsumura-Shi-1]},
\cite{[Huang-Qin]}, \cite{[Kawashima-Nishibata-Zhu]},
\cite{[Matsumura]} \cite{[Matsumura-Nishihara-3]}, \cite{[Qin]},
\cite{[Zhu]}, etc.. For the inflow problem  of the full Naier-Stokes
equation (\ref{(1.1)})--(\ref{(1.4)}), Huang-Li-Shi
\cite{[Huang-Li-Shi]} proved the stability of the BL-solution in
some cases. More precisely, they show that when
$(\r_\pm,u_\pm,\t_\pm)$ both belong to the subsonic region, the
BL-solution is expected and the stability of this BL-solution and
its superposition with the 3-rarefaction wave is proved under some
smallness assumptions. Notice that both the BL-solution and the
rarefaction wave must be weak enough. When the boundary value
$(\r_-,u_-,\t_-)$ belongs to the supersonic region, there is no
BL-solution. Thus the large-time behavior of the solution is
expected to be same as that of the Cauchy problem and the stability
of the rarefaction waves is given.

In this paper, firstly we give the existence (or non-existence) of
the BL-solution when the right end state $(\r_+,u_+,\t_+)$ belongs
to the subsonic, transonic and supersonic regions, respectively. The
rigorous proof is given in Appendix. Notice that it is more natural
to present the classifications according to the locations of the
right end state $(\r_+,u_+,\t_+)$ from the qualitative theory of the
autonomous ODE system. Then we prove the stability of not only the
single contact wave but also the superposition of the BL-solution
(subsonic case), the viscous contact wave and the 3-rarefaction wave
to the inflow problem (\ref{(1.1)})--(\ref{(1.4)}). Here the
amplitude of the rarefaction wave can be arbitrarily large.

Now we briefly review some key analytic techniques in studying the
stability of the basic wave patterns. The strict monotonicity of the
corresponding characteristic speed along the wave profiles plays a
crucial role in stability analysis of the viscous shock wave and
rarefaction wave. Precisely speaking, the shock wave is a
compression wave, so the characteristic speed is monotone decreasing
in the shock profile. Thus anti-derivative variable to the
perturbation should be introduced in the stability analysis. While
the rarefaction wave is an expansion wave and the characteristic
speed is monotone increasing along the rarefaction wave, thus the
direct energy estimates to the perturbation itself are available.
However, the characteristic speed along the contact wave is
constant, and the spatial derivative of the velocity changes signs
along the contact wave profile. Due to the degenerate
characteristics, the stability of the contact wave profile to the
compressible Navier-Stokes system (\ref{(1.1)}) is just proved by
\cite{[Huang-Matsumura-Xin-1]} and \cite{[Huang-Xin-Yang]} in 2005,
twenty years later than the nonlinear wave  in 1985.  In
\cite{[Huang-Matsumura-Xin-1]} and \cite{[Huang-Xin-Yang]}, the
anti-derivative variable to the perturbation is introduced and the
proof framework is motivated by the viscous shock profile. Notice
that a convergence rate of the order of $(1+t)^{-\f14}$ in sup-norm
is a by-product of the estimation. However, there is no convergence
rate obtained so far for the viscous shock wave and the rarefaction
wave.

Recently, Huang-Matsumura-Xin \cite{[Huang-Matsumura-Xin-2]}
obtained a new estimate on the heat kernel which can be applied to
the study of the stability of the viscous contact wave in the
framework of the rarefaction wave, see \cite{[Huang-Li-Matsumura]}
or Lemma 3.4 in the present paper. Namely, the anti-derivative
variable of the perturbation is not needed and the estimations to
the perturbation itself are also suitable to get the stability of
the viscous contact wave. But the time-decay rate can not be gotten
as a compensation. More importantly, the advantage of this framework
is that it can be used to study the stability of the contact wave to
the IBVP of (\ref{(1.1)}) since the boundary terms could be treated
conveniently. We will make full use of this new estimate on heat
kernel to study the inflow problem $(\ref{(1.1)})$ and get the
expected results.

The novelty of the paper lies in the following three aspects: (1)
The rigorous proof and the classifications of the existence (or
non-existence) of the BL-solution to the inflow problem. (2) The
stability of the superposition of three different wave patterns (the
BL-solution, the viscous contact wave and the rarefaction wave). (3)
The large amplitude of the rarefaction wave in the superposition
wave. The main difficulties in our proofs are how to deal with the
boundary terms and the interactions of three different wave
patterns.

Because the system $(\ref{(1.1)})$ we consider  is in one dimension
of the space variable $x$, it is convenient to use the following
Lagrangian coordinate transformation:
$$
x\Rightarrow \int_0^x \r (y,t)dy,\qquad t\Rightarrow t.
$$
Thus the system $(\ref{(1.1)})$ can be transformed into the
following moving boundary problem of Navier-Stokes equations in the
Lagrangian coordinates:
\begin{equation}
\left\{
\begin{array}{ll}
\di v_t-u_x=0,&x>\s_-t,~~t>0,\\[1mm]
\di u_t+p_x=\m (\f{u_x}{v})_x,&x>\s_-t,~~t>0,\\[1mm]
\di
(e+\f{u^2}{2})_t+(pu)_x=\k(\f{\t_x}{v})_x+\m(\f{uu_x}{v})_x,~&x>\s_-t,~~t>0,\\[3mm]
\di (v,u,\t)(t,x=\s_-t)=(v_-,u_-,\t_-),&u_->0,\\[2mm]
\di (v,u,\t)(t=0,x)=(v_0,u_0,\t_0)(x)\ra(v_+,u_+,\t_+),&{\rm as}~~
x\ra+\i,
\end{array}
\right. \label{(1.5)}
\end{equation}
 where $ v(t,x)=\f{1}{\r(t,x)}$ represents the specific
volume of the gas, and the   boundary moves with the constant speed
$\s_-=-\f{u_-}{v_-}<0$. Now we have that for the perfect gas,
\begin{equation}
p=\f {R\t} v=Av^{-\g}\exp{(\f{\g-1}{R}s)},\qquad e=\f{R}{\g-1}\t
+{\rm const.}\label{(1.6)}
\end{equation}
In order to fix the moving boundary $x=\s_-t$, we introduce a new
variable $\x=x-\s_-t$. Then we have the half space problem
\begin{equation}
\left\{
\begin{array}{ll}
\di v_t-\s_-v_\x-u_\x=0,&\x>0,~~t>0,\\[1mm]
\di u_t-\s_-u_\x+p_\x=\m (\f{u_\x}{v})_\x,&\x>0,~~t>0,\\[2mm]
\di
(e+\f{u^2}{2})_t-\s_-(e+\f{u^2}{2})_\x+(pu)_\x=\k(\f{\t_\x}{v})_\x+\m(\f{uu_\x}{v})_\x,~&\x>0,~~t>0,\\[3mm]
\di (v,u,\t)(t,\x=0)=(v_-,u_-,\t_-),&u_->0,\\[2mm]
\di (v,u,\t)(t=0,\x)=(v_0,u_0,\t_0)(\x)\ra(v_+,u_+,\t_+),&{\rm
as}~~\x\ra+\i.
\end{array}
\right. \label{(1.5+)}
\end{equation}
Given the right end state $(v_+,u_+,\t_+)$, we can define the
following wave curves in the phase space $(v,u,\t)$ with $v>0$ and
$\t>0$.

$\bullet$ Contact wave curve:
\begin{equation}
CD(v_+,u_+,\t_+)= \{(v,u,\t)  |  u=u_+, p=p_+, v \not\equiv v_+
 \}. \label{(1.7)}
\end{equation}

$\bullet$ BL-solution curve (subsonic case, i.e., $(v_+,u_+,\t_+)\in
\Omega_{sub}^+$):
\begin{equation}
BL(v_+,u_+,\t_+)=\bigg\{(v,u,\t)\bigg{ |}
\f{u}{v}=-\s_-=\f{u_+}{v_+}, (u,\t)\in \mathcal{M}(u_+,\t_+)
\bigg\},\label{(1.8)}
\end{equation}
where $\mathcal{M}(u_+,\t_+)$ is the center-stable manifold defined
in Lemma 2.1 below.

$\bullet$ 3-Rarefaction wave curve:
\begin{equation}
 R_3 (v_+, u_+, \theta_+):=\Bigg{ \{} (v, u, \theta)\Bigg{ |}v>v_+ ,~u=u_+-\int^v_{v_+}
 \lambda_3(\eta,
s_+) \,d\eta,~ s(v, \theta)=s_+\Bigg{ \}},\label{(1.9)}
\end{equation}
where $s_+=s(v_+,\t_+)$ and $\l_3=\l_3(v,s)$ is the third
characteristic speed given in (\ref{(2.1)}).

\vskip 1mm
 Our main stability results are, roughly speaking, as
follows:

\vskip 1mm
 (I). If the state $(v_-,u_-,\t_-) \in {\rm CD}
(v_+,u_+,\t_+)$, then the viscous contact wave is asymptotic stable
under some smallness conditions which are given in Theorem 2.1.

\vskip 1mm
 (II). If the state $(v_-,u_-,\t_-) \in {\rm
BL\texttt{-}CD\texttt{-}R_3} (v_+,u_+,\t_+)$, then there exist a
unique state  $(v_*,u_*,\t_*)\in\Omega_{sub}^+$ and a unique state
$(v^*,u^*,\t^*)$, such that $(v_-,u_-,\t_-) \in {\rm
BL}(v_*,u_*,\t_*)$, $(v_*,u_*,\t_*) \in {\rm CD}(v^*,u^*,\t^*)$ and
$(v^*,u^*,\t^*)\in {\rm R_3} (v_+,u_+,\t_+)$ and the superposition
of the BL-solution, the viscous contact wave and the rarefaction
wave is asymptotically stable provided that $|(u_--u_*,\t_--\t_*)|$
and $|v_*-v^*|$ are suitably small and the conditions in Theorem 2.2
hold. It is remarked that the BL-solution and the viscous contact
wave must be weak but the rarefaction wave is not necessarily weak.

\vskip 1mm
 {\bf Notations:} Throughout the paper several positive
generic constants are denoted by $c,c_0,c_1,\cdots$ or
$C,C_1,C_2,\cdots$ without confusions. The small constant $\nu>0$ is
used in Young inequality
$$
ab\leq \nu a^{p_1}+C_\nu b^{p_2},\qquad \f{1}{p_1}+\f{1}{p_2}=1,
$$
where $C_\nu$ is the constant depending on $\nu$. For functional
space, $H^l(\mathbf{R}^+)$ denotes the $l$-order Sobolev space with
the norm
$$
\|f\|_l=\left(\sum_{i=0}^l\|\partial_x^i f\|^2\right)^{\f12},\quad
{\rm where} \|f\|=\|f\|_{L^2}.
$$
\section{Preliminaries and main results}
\setcounter{equation}{0} In this section, we will show the wave
patterns considered in the paper. We start with the BL-solution.

\subsection{BL-solution}
It is known that the hyperbolic system (\ref{(1.4+)}) has three
characteristic speeds
\begin{eqnarray}
\lambda_1= -\sqrt{\f{\g p}{v}} ,~~~ \lambda_2=0 ,~~~ \lambda_3=
\sqrt{\f{\g p}{v}}. \label{(2.1)}
\end{eqnarray}
The sound speed $c(v,\t)$ and the Mach number $M$ are defined by
\begin{equation}
c(v,\t)=v\sqrt{\f{\g p}{v}}=\sqrt{R\g \t}, \label{(2.2)}
\end{equation}
and
\begin{equation}
M(v,u,\t)=\f{|u|}{c(v,\t)}=\f{|u|}{\sqrt{R\g\t}},\label{(2.3)}
\end{equation}
respectively.

We divide the phase space $\{(v,u,\theta), v>0, \theta>0\}$ into
three regions:
\begin{equation}
\begin{array}{ll}
\Omega_{sub}&\di :=\big{\{}(v, u,\theta)|~M(v,u,\t)<1~\big{\}}, \\
\Gamma_{trans}&\di :=\big{\{}(v, u,\theta)| ~M(v,u,\t)=1~\big{\}}, \\
\Omega_{super}&\di:= \big{\{}(v,u,\theta)|~ M(v,u,\t)>1 ~\big{\}}.
\end{array}
\label{(2.4)}
\end{equation}
Call them the subsonic, transonic and supersonic regions,
respectively. If adding the alternative condition $u>0$ or $u<0$,
then we have six connected subsets $\Omega_{sub}^{\pm}$,
$\Gamma_{trans}^{\pm}$, and $\Omega_{super}^{\pm}$.

When $(v_-,u_-,\t_-)\in \Omega_{sub}^+$, we have
$$
\l_1(v_-,u_-,\t_-)<\s_-<0,
$$
hence the existence of the traveling wave solution
\begin{equation}
\begin{array}{l}
(V^B,U^B,\T^B)(\x),\quad \x=x-\s_-t,\\[2mm]
(V^B,U^B,\T^B)(0)=(v_-,u_-,\t_-),\quad
(V^B,U^B,\T^B)(+\i)=(v_+,u_+,\t_+),
\end{array} \label{(2.5)}
\end{equation}
to (\ref{(1.5)}), or the stationary solution to (\ref{(1.5+)}) is
expected. We call this traveling wave solution $(V^B,U^B,\T^B)(\x)$
the boundary layer solution to the inflow problem (\ref{(1.5)}).
Note that the speed of the traveling wave is just the speed of the
moving boundary of (\ref{(1.5)}).

 In the following, we will
give the existence (or non-existence) of BL-solution to the inflow
problem (\ref{(1.5)}). From (\ref{(2.5)}), BL-solution $(V^B, U^B,
\Theta^B)(\x)$ satisfies the following ODE system
\begin{equation}
\left\{
\begin{array}{ll}
\di -\s_-(V^{B})^\prime-(U^{B})^\prime=0 , \qquad\qquad\qquad ^\prime:=\frac{d}{d\x}&\x>0, \\
\di -\s_-(U^{B})^\prime+(P^{B})^\prime= \mu\big{(}\frac{(U^{B })^\prime}{V^B}\big{)}^\prime, &\x>0, \\
\di   -\s _-(\f{R}{\g-1}\Theta^B+\frac{(U^B)^2}{2})^\prime+
(P^BU^B)^\prime =\k( \f{(\T^{B})^\prime}{V^B})^\prime +\m (\f{U^B(U^
{B})^\prime
 }{V^B})^\prime  ,& \x>0 , \\
\di  (V^B, U^B, \Theta^B)(0)=(v_-, u_-, \t_-) ,\quad (V^B, U^B,
\Theta^B)(+\infty)=(v_+, u_+, \t_+).
\end{array}
\right. \label{(2.6)}
\end{equation}
where $P^B:=p(V^B, \Theta^B)=\f{R\T^B}{V^B}$.

Integrating the system (\ref{(2.6)}) over $(\xi,+\i)$ implies that
\begin{equation}
\left\{
\begin{array}{ll}
\di  -\s_-(V^B-v_+)-(U^B-u_+)=0,\\
\di \m\frac{(U^{B})^\prime}{V^B}=-\s_-(U^{B}-u_+)+R\Bigg{(} \frac{\T^{B}}{V^B}- \frac{\t_+}{v_+}\Bigg{)}   , \\
\di   \k \f{(\T^{B})^\prime}{V^B}=-\s _- \f{R}{\g-1}(\Theta^B-\t_+)+
 p_+(U^B-u_+)+\frac{\s_-}{2}(U^B-u_+)^2  , \\
\di  (U^B, \Theta^B)(0)=(  u_-, \t_-) ,\quad (  U^B,
\Theta^B)(+\infty)=(  u_+, \t_+).
\end{array}
\right. \label{(2.7)}
\end{equation}
Let $\x=0$ in $(\ref{(2.7)})_1$, we have
\begin{equation}
\s_-=-\f{u_-}{v_-}=-\f{U^B}{V^B}=-\f{u_+}{v_+},\label{(2.8)}
\end{equation}
which is the first condition of BL-solution curve in (\ref{(1.8)}).

From the fact $u_->0$ and $v_{\pm}>0$, we find that $u_+$ must
satisfy
\begin{equation}
u_+=\f{u_-}{v_-} v_+>0.\label{(2.9)}
\end{equation}
Rewrite $(\ref{(2.7)})_2$--$(\ref{(2.7)})_4$ as
\begin{equation} \left\{
\begin{array}{ll}
\di  (U^{B })^\prime=-\frac{\s_-}{\m}V^B(U^{B}-u_+)+\frac{R}{\m}\bigg{(}  \T^{B}  - \frac{\t_+}{v_+}V^B\bigg{)}   , \\
\di  (\T^{B })^\prime =- \f{R\s_-V^B}{\k(\g-1)}(\Theta^B-\t_+)+
 \frac{p_+}{\k}V^B(U^B-u_+)+\frac{\s_-V^B}{2\k}(U^B-u_+)^2  , \\
\di  (  U^B, \Theta^B)(0)=( u_-, \t_-) ,\quad (  U^B,
\Theta^B)(+\infty)=(  u_+, \t_+).
\end{array}
\right. \label{(2.10)}
\end{equation}
Denote
\begin{eqnarray}
 \bar{U}^B =U^B-u_+, \quad   \bar{\Theta}^B =   \Theta^B -\theta_+ ,\label{(2.11)}
\end{eqnarray}
and
\begin{equation}
J=\left ( \begin{array}{cc} \frac{u_+^2-R\t_+}{\m u_+}& \frac{R}{
\m}\\ [3mm]\frac{R\t_+}{\k} &\f{Ru_+}{\k(\g-1)}
\end{array}
\right)=\left ( \begin{array}{cc}
  \frac{(M_+^2\gamma-1)u_+}{M_+^2\g \m}
&\frac{R}{\m}\\[3mm]  \frac{u_+^2 }{ M_+^2  \gamma \k} &\frac{R u_+ }{
\k(\g-1)}
\end{array}
\right) ,\label{(2.12)}
\end{equation}
where $M_+=M(v_+,u_+,\t_+)$.

Then we obtain the automatous ODE system
\begin{eqnarray}
\left\{
\begin{array}{l}
\left(
\begin{array}{c}
\bar{U}^B\\ \bar\Theta^B
\end{array}
\right)^\prime=J\left(
\begin{array}{c}
\bar{U}^B\\ \bar\Theta^B
\end{array}
\right)+\left(
\begin{array}{c}
F_1(\bar{U}^B,\bar\Theta^B)\\
F_2(\bar{U}^B,\bar\Theta^B)
\end{array}
\right)\\[5mm]
 (\bar U^B, \bar\Theta^B)(0)=(u_--u_+, \t_--\t_+) ,\quad (\bar U^B,
\bar\Theta^B)(+\infty)=( 0, 0).
\end{array}
\right.
 \label{(2.13)}
\end{eqnarray}
where
\begin{equation}
\begin{array}{ll}
\di F_1(\bar{U}^B,\bar\Theta^B)&= \di \frac{1 }{\m}  (\bar{U}^B)^2, \\
\di F_2(\bar{U}^B,\bar\Theta^B)&= \di \Bigg{(}
 \frac{R\t_+}{\k u_+}-\frac{u_+}{2\k}\Bigg{)} (\bar{U}^B)^2+\f{R  }{\k(\g-1)} \bar{U}^B \bar{\Theta}^B -\frac{ 1}{2\k}
(\bar{U}^B)^3.
\end{array}
\label{(2.14)}
\end{equation}
Now we state the existence results of the solution to (\ref{(2.13)})
while its   proof will be shown in Appendix.

\vskip 2mm

\noindent\textbf{Lemma 2.1 (Existence of BL-solution)} \emph{
Suppose that $v_\pm>0$, $u_- >0$, $\theta_\pm>0$ and let
$\d^B=|(u_+-u_-,\t_+-\t_-)|$. If $u_+\leq 0$, then there is no
solution to (\ref{(2.10)}) or (\ref{(2.13)}).  If $u_+>0$, then
there exists a suitably small constant $\delta>0$
 such that if $0<\d^B\leq\delta$,   then}

\emph{Case I. Supersonic case: $M_+>1$. Then   there is no solution
to (2.10) or (2.13).}

\emph{Case II.  Transonic case: $M_+=1$. Then
 $(u_{+}, \theta_+)$ is a saddle-knot point to
 (\ref{(2.13)}). Precisely,  there exists a unique trajectory $ \Gamma $ tangent to the line
 $$
 \m u_+(U^B-u_+)- \k(\g-1)(\Theta^B-\theta_+)=0
$$
 at the point $(u_{+}, \theta_+)$. For each $( u_-,\theta_-)\in \Gamma $, there exists a unique solution
 $(U^B, \Theta^B) $ satisfying}
\begin{eqnarray}
\Bigg{|}\f{d^{ n }}{d\x^n}(U^B-  u_+ ,  \Theta^B - \theta_+
)\Bigg{|} \leq C\frac{(\d^B)^n}{(1+\d^B \x)^n} ,~~~n=0,1, 2,
\dots,~~ ~~\x\in \mathbb{R}_+ .\label{(2.15)}
\end{eqnarray}

\emph{Case III.  Subsonic case: $M_+<1$. Then the equilibrium point
$(u_{+},\theta_+)$ is a saddle point of (\ref{(2.13)}). Precisely,
there exists a center-stable manifold $\mathcal{M} $ tangent to the
line
$$
(1+a_2c_2u_+)(U^B-u_+) -a_2(\Theta^B-\theta_+)=0
$$
on the opposite directions at the point $(u_{+}, \theta_+)$. Here
$c_2$ is one of the solutions to the equation
$$
 y^2+ \Bigg{(}\frac{M_+^2\gamma-1}{M_+^2R\gamma}-\frac{\m }{\k(\g-1)}\Bigg{)} y-\frac{\m }{M_+^2 R \gamma \k}=0
$$
 and $a_2=-\frac{R}{\m(\l_A^1-\l_A^2) }$ with $\l_A^1>0,~~\l_A^2<0$ are the two eigenvalues of the matrix $A$.
 Only when $(u_-, \theta_-)\in \mathcal{M} $, does there
 exist  a unique solution $(U^B, \Theta^B )\subset \mathcal{M} $   satisfying}
\begin{eqnarray}
\Bigg{|}\f{d^{ n }}{d\x^n}(U^B - u_+, \Theta^B - \theta_+ )\Bigg{|}
\leq C \d^B  e^{-c\x }  ,~~~n=0,1, 2, \dots,~~~~~\x\in
\mathbb{R}_+.\label{(2.16)}
\end{eqnarray}

{\bf Remark:} This Lemma is the first one for the classifications of
the BL-solution to the inflow problem (1.1). The stability of the
single BL-solution in subsonic case (Case III) is proved in
\cite{[Huang-Li-Shi]}. In this paper, we are concerned with the
stability of the superposition of this subsonic BL-solution with
viscous contact wave and rarefaction wave. As for the stability of
the BL-solution in transonic case (Case II) and its superposition
with other wave patterns, it is in consideration \cite{[Qin-Wang]}.

\subsection{Viscous contact wave}

If $(v_-,u_-,\t_-)\in CD(v_+,u_+,\t_+)$, i.e.,
\begin{equation}
u_-=u_+,~p_-=p_+,\label{(2.17)}
\end{equation}
then the following Riemann problem of the Euler system
\begin{equation}
\left\{
\begin{array}{ll}
v_t-u_x=0,\\
u_t+p_x=0,\\
(e+\f{u^2}{2})_t+(pu)_x=0,
\end{array}
\right.\label{(2.18)}
\end{equation}
with Riemann initial value
$$
(v,u,\t)(t=0,x)=\left\{
\begin{array}{ll}
(v_-,u_-,\t_-),\qquad & x<0,\\
(v_+,u_+,\t_+),\qquad & x>0,
\end{array}
\right.
$$
admits a single contact discontinuity solution
$$
(v,u,\t)(t,x)=\left\{
\begin{array}{ll}
(v_-,u_-,\t_-),\qquad & x<0,~ t>0,\\
(v_+,u_+,\t_+),\qquad & x>0,~ t>0.
\end{array}
\right.
$$

From \cite{[Huang-Matsumura-Xin-1]}, the viscous version of the
above contact discontinuity, called viscous contact wave
$(V^{CD},U^{CD},\T^{CD})(t,x)$, could be defined by
\begin{equation}
\begin{array}{ll}
\di \T^{CD}(t,x)=\T^{Sim}(\f{x}{\sqrt{1+t}}),\\
\di V^{CD}(t,x)=\f{R\T^{CD}(t,x)}{p_+},\\
\di U^{CD}(t,x)=u_+
+\f{\k(\g-1)}{R\g}\f{\T^{CD}_{x}(t,x)}{\T^{CD}(t,x)},
\end{array}
\label{(2.19)}
\end{equation}
where $\T^{Sim}(\eta)$, $\di\eta=\f{x}{\sqrt{1+t}}$, is the unique
self-similar solution of the nonlinear diffusion equation
\begin{equation}
\left\{
\begin{array}{l}
\di \T_t=\f{\k p_+(\g-1)}{R^2\g}\left(\f{\T_x}{\T}\right)_x,\\
\di \T(\pm\i)=\t_\pm.
\end{array}
\right. \label{(2.20)}
\end{equation}
Thus the viscous contact wave defined in (\ref{(2.19)}) satisfies
the following property
\begin{equation}
|\T^{CD}-\t_\pm|+(1+t)^{\f12}|\T^{CD}_x|+(1+t)|\T^{CD}_{xx}|+(1+t)^{\f32}
|\T^{CD}_{xxx}| =O(1)\d^{CD} e^{-\f{c_0x^2}{1+t}}, \label{(2.21)}
\end{equation}
as $|x|\ra+\i$, where $\d^{CD}=|\t_+-\t_-|$ is the amplitude of the
viscous contact wave and $c_0$ is a positive constant. Note that
$\x=x-\s_-t$, then
$(V^{CD},U^{CD},\T^{CD})(t,x)=(V^{CD},U^{CD},\T^{CD})(t,\x+s_-t)$
satisfies the system
\begin{equation}
\left\{
\begin{array}{l}
\di V^{CD}_{t}-\s_-V^{CD}_{\x}-U^{CD}_{\x}=0,\\
\di U^{CD}_{t}-\s_-U^{CD}_{\x}+P^{CD}_{\x}=\m(\f{U^{CD}_{\x}}{V_{C}})_\x+\bar Q_1,\\
\di
\f{R}{\g-1}(\T^{CD}_{t}-\s_-\T^{CD}_{\x})+P^{CD}U^{CD}_{\x}=\k(\f{\T^{CD}_{\x}}{V^{CD}})_\x+\m\f{(U^{CD}_{\x})^2}{V^{CD}}+\bar
Q_2,
\end{array}
\right. \label{(2.22)}
\end{equation}
where $\di P^{CD}=\f{R\T^{CD}}{V^{CD}}=p_+=p_-$ and the error terms
$\bar Q_1,~\bar Q_2$ are given by
\begin{equation}
\begin{array}{ll}
\di \bar Q_1&\di =(U^{CD}_{t}-\s_-U^{CD}_{\x})-\m
(\f{U^{CD}_{\x}}{V^{CD}})_\x=O(1)(|\T^{CD}_{\x}|^3+|\T^{CD}_{\x\x\x}|+|\T^{CD}_{\x\x}||\T^{CD}_{\x}|)\\
\di &\di =O(1)\d^{CD}
(1+t)^{-\f32}e^{-\f{c_0(\x+\s_-t)^2}{1+t}},\qquad {\rm
as}~~|\x+\s_-t|\ra +\i,\\[4mm]
 \di \bar Q_2&\di =-\m
\f{(U^{CD}_{\x})^2}{V^{CD}}=O(1)(|\T^{CD}_{\x}|^4+|\T^{CD}_{\x\x}|^2)\\
\di &\di =O(1)(\d^{CD})^2
(1+t)^{-2}e^{-\f{c_0(\x+\s_-t)^2}{1+t}},\qquad {\rm
as}~~|\x+\s_-t|\ra +\i.
\end{array}
\label{(2.23)}
\end{equation}

\subsection{Rarefaction wave}

If $(v_-, u_-, \theta_-) \in R_3 (v_+, u_+, \theta_+)$, then there
exists  a 3-rarefaction wave $(v^r, u^r, s^r)(x/t)$ which is the
global (in time) weak solution of the following Riemann problem
\begin{eqnarray}\label{4}
\left\{
\begin{array}{l}
\di  v_{t}^r- u_{x}^r= 0,\\
\di u_{t}^r +  p_{x}(v^r, \theta^r) = 0 , \qquad
\qquad\qquad\qquad~~\,   t>0, x\in \textbf{R},
\\
\di \frac{R}{\g-1}\theta_t^r + p(v^r, \theta^r)
u_x^r =0,\\
\di  (v^r, u^r, \t^r)(0,x)=\left\{
\begin{array}{l}
\di (v_-, u_-, \t_-),   x<
0 ,\\
\di  (v_+, u_+, \t_+), x> 0 .
\end{array}
\right.
\end{array} \right.\label{(2.24)}
\end{eqnarray}

From \cite{[Huang-Matsumura-Shi-3]}, it is convenient to construct
the approximated rarefaction wave $(V^R, U^R, \Theta^R)(t,x)$ to the
inflow problem (\ref{(1.5)}) by the solution of the Burgers equation
\begin{eqnarray}\label{4}
\left\{
\begin{array}{l}
\di w_{t}+ww_{x}=0,\\
\di w( 0,x )=w_0(x)=\left\{
\begin{array}{l}
\di w_-,\qquad \qquad\qquad\qquad~~\,  x<
0 ,\\
\di   w_-+C_q \d^r \int^{ \varepsilon x }_0y^qe^{-y}\,dy, x\geq 0 ,
\end{array}
\right.
\end{array}
\right.\label{(2.25)}
\end{eqnarray}
where $\d^r=w_+-w_-$, $q\geq 16$ is some fixed constant, $C_q$ is a
constant such that $C_q\int^{\infty}_0y^qe^{-y} dy=1 $, and
$\varepsilon\leq 1$ is a small positive constant to be determined
later. The solution of the Burgers equation (\ref{(2.25)}) have the
following properties:
 \vskip 2mm
\noindent\textbf{Lemma 2.2 (\cite{[Huang-Matsumura-Shi-3]})}
\emph{Let $0<w_-<w_+$, Burgers equation $(\ref{(2.25)})$ has a
unique smooth solution $w(t,x)$ satisfying}

i) ~ $w_-\leq w(t,x)<w_+,~w_x(t,x)\geq 0 $,

ii) ~ For any $p$ $(1\leq p\leq \infty)$, there exists a constant
$C_{pq}$ such that
\begin{eqnarray}\nonumber
&& \parallel w_x(t)\parallel_{L^p}\leq
C_{pq}\min\big{\{}\delta_r\varepsilon^{1- 1/p},~
\delta_r^{1/p}t^{-1+1/p}\big{\}},   \cr&&\parallel
w_{xx}(t)\parallel_{L^p}\leq
C_{pq}\min\big{\{}\delta_r\varepsilon^{2-1/p},~ \big{(}
\delta_r^{1/p }+\delta_r^{ 1/q}\big{)}t^{-1+1/q}\big{\}},
\end{eqnarray}
\indent iii) ~ If $ x< w_-t $, then $w(t,x)\equiv w_-$,

iv) ~ $\sup\limits_{x\in\mathbb{R}} |w(t,x)-w^r(x/t)|\rightarrow 0$,
as $t\rightarrow \infty$.

Thus we construct the approximated rarefaction wave $(V^R, U^R,
\Theta^R) (t,x)$ by
\begin{eqnarray}\label{4}
\left\{
\begin{array}{l}
\di  S^R=s(V^R,\T^R)=s_+,\\
\di w(1+t,x)= \l_3(V^R(t,x),s_+),
\\\di U^R(t,x)=u_+-\int^{V^R(t,x)}_{v_+}  \l_3(v,s_+) dv.
\end{array} \right.\label{(2.26)}
\end{eqnarray}
Note that $\x=x-\s_-t$, then the smoothed 3-rarefaction wave $(V^R,
U^R, \Theta^R)(t,\x)$ defined above satisfies
\begin{equation}
\left\{
\begin{array}{l}
\di V^R_{t}-\s_-V^R_{\xi}-U^R _{\xi}=0,\\
\di U^R_{t}-\s_-U^R_{\xi}+P^R_{\xi}= 0,  \qquad
\qquad\qquad   \x>0,~t>0, \\
\di \f{R}{\g-1}(\T^R _{t}-\s_-\Theta^R_{\xi})+P^RU^R_{\xi}=0,\\
\di  (V^R, U^R, \Theta^R)(t,0)=(v_-, u_-, \t_-),\quad (V^R, U^R,
\Theta^R)(t,+\infty)=(v_+, u_+, \t_+),
\end{array}
\right.\label{(2.27)}
\end{equation}
where $P^R:=p(V^R, \Theta^R)$.

\vskip 2mm

\noindent\textbf{Lemma 2.3 (\cite{[Huang-Matsumura-Shi-3]})} Let
$\delta^R=|(v_+, u_+, \t_+)-(v_-, u_-, \t_-)|$. The approximated
rarefaction wave $(V^R, U^R, \T^R)(t,\x)$ satisfies

  i) For $\x>0,~,t>0,$ $ U^R_{\x}(t,\x)\geq 0 $,

   ii) For any $p$  ($1\leq p\leq \infty$),  there exists a constant
$C_{pq}$ such that for $t\geq 0$,
\begin{eqnarray}\nonumber
&& \parallel (V^R_{\x}, U^R_{\x}, \T^R_{\x})(t)\parallel_{L^p}\leq
C_{pq}\min\big{\{}\delta^R \varepsilon^{1- 1/p },~ (\delta^R)^{1/p
}(1+t)^{-1+1/p}\big{\}},  \cr&&\parallel (V^R_{\x\x},U^R_{\x\x},
\T^R_{\x\x})(t)\parallel_{L^p}\leq C_{pq}\min\big{\{}  \delta^R
\varepsilon^{2-1/p},~ \big{(} (\delta^R)^{1/p }+(\delta^R)^{
1/q}\big{)}(1+t)^{-1+ 1 /q }\big{\}},
 \end{eqnarray}

  iii) If $ \x+\s_- t\leq \l_3(v_-, u_-,
  \t_-)(1+t) $, then $(V^R, U^R, \T^R)(t,\x)\equiv (v_-, u_-,
  \t_-)$,

   iv) $ \sup\limits_{\x\in\mathbb{R}_+}\big{|}(V^R, U^R, \T^R)( t, \x)-(v^r,u^r,\t^r)
\big{(}\frac{\x}{1+t}\big{)}\big{|}\rightarrow 0$, as $t\rightarrow
\infty$.

\subsection{Main results}
Now we can state our main results. The first one is in the
following.

\vskip 2mm

\noindent{\bf Theorem 2.1} If $(v_-,u_-,\t_-)\in CD(v_+,u_+,\t_+)$.
Let $(V^{CD},U^{CD},\T^{CD})(t,x)$ be the viscous contact wave
defined in (\ref{(2.19)}). There exists a small constant $\d_0$ such
that if the wave amplitude $\d^{CD}$ and the initial values satisfy
$$
\d^{CD}+\|(v_0-V^{CD}_{0},u_0-U^{CD}_{0},\t_0-\T^{CD}_{0})\|_{1}\leq
\d_0,
$$
then the moving boundary problem (\ref{(1.5)}) or the half space
problem (\ref{(1.5+)}) admits a unique global solution
$(v,u,\t)(t,\x)$ satisfying
$$
\begin{array}{l}
\di (v-V^{CD},u-U^{CD},\t-\T^{CD})(t,\x)\in C([0,+\i),H^1({\bf R}^+)),\\[2mm]
\di (v-V^{CD})_\x(t,\x)\in L^2(0,+\i;L^2({\bf R}^+)),\\[2mm]
\di ((u-U^{CD})_\x,(\t-\T^{CD})_\x)(t,\x)\in L^2(0,+\i;H^1({\bf
R}^+)),
\end{array}
$$
and
$$
\lim_{t\rightarrow +\i}\sup_{\x\in {\bf
R}^+}|(v-V^{CD},u-U^{CD},\t-\T^{CD})(t,\x)|=0.
$$

Now we state our second result. If $(v_-,u_-,\t_-) \in {\rm
BL\texttt{-}CD\texttt{-}R_3} (v_+,u_+,\t_+)$, then there exist
states $(v_*,u_*,\t_*)\in\Omega_{sub}^+$ and $(v^*,u^*,\t^*)$, such
that $(v_-,u_-,\t_-) \in {\rm BL}(v_*,u_*,\t_*)$, $(v_*,u_*,\t_*)
\in {\rm CD}(v^*,u^*,\t^*)$ and $(v^*,u^*,\t^*)\in {\rm R_3}
(v_+,u_+,\t_+)$. In fact, by$(v_-,u_-,\t_-) \in {\rm
BL}(v_*,u_*,\t_*)$ and (\ref{(1.8)}), we have
\begin{equation}
u_*=-\s_-v_*,\label{(2.52)}
\end{equation}
and by $(v^*,u^*,\t^*)\in {\rm R_3} (v_+,u_+,\t_+)$, (\ref{(1.9)})
gives
\begin{equation}
u^*=u_+-\int_{v_+}^{v^*}\l_3(\eta,s_+)d\eta,\quad v^*>v_+.
\label{(2.53)}
\end{equation}
Thus the two curves (\ref{(2.52)}) and (\ref{(2.53)}) have a unique
intersection point $u=\tilde{u}$ in $(v,u)$ space. If
$u_*=u^*=\tilde{u}$, then $v_*=v^*$, thus there is no contact wave.
By $(v_*,u_*,\t_*) \in {\rm CD}(v^*,u^*,\t^*)$, we have
$$
u_*=u^*,\qquad v_*\ne v^*.
$$
Thus if $u_*=u^*\ne \tilde{u}$, then $v_*\ne v^*$ and there exists a
contact wave. Among the three values $u_*(=u^*)$, $v_*$ and $v^*$,
only one is independent, the other two can be determined
accordingly.

Now assume that $u_*(=u^*)$ is given, then from (\ref{(2.52)}) and
(\ref{(2.53)}), we can determine $v_*$ and $v^*$ by
$$
v_*=\f{u_*}{-\s_-},\qquad
v^*=v_+\left[\f{\g-1}{2A\sqrt{R\g\t_+}}(u_*-u_+)+1\right]^{\f{2}{1-\g}}.
$$
By the definition of the rarefaction wave curve $R_3$ in
(\ref{(1.9)}), we have $s(v_*,\t_*)=s_+$, i.e.,
$$
\t^*=\t_+\left(\f{v_+}{v^*}\right)^{\g-1}=\t_+\left[\f{\g-1}{2A\sqrt{R\g\t_+}}(u_*-u_+)+1\right]^2.
$$
Again by the contact wave curve (\ref{(1.7)}), one must have
$$
p_*=p^*,
$$
thus
\begin{equation}
\t_*=\f{\t^*}{v^*}v_*=\f{\t_+}{-\s_-v_+}u_*\left[\f{\g-1}{2A\sqrt{R\g\t_+}}(u_*-u_+)
+1\right]^{\f{2\g}{\g-1}}.\label{(2.30+)}
\end{equation}
So if $u_*$ large enough, then $(u_*,\t_*)$ in (\ref{(2.30+)}) must
belong to the region $\Omega_{sub}^+:=\{(u_*,\t_*)|0<u_*<\sqrt{R\g
\t_*}\,\}$. Moreover, from the definition of the BL-solution,
$(u_-,\t_-)\in \mathcal{M}(u_*,\t_*)$. Thus $(u_*,\t_*)$ can be
determined uniquely if $(u_-,\t_-)$ is given suitably.

 Define the superposition wave $(V,
U, \Theta)(t,\x)$ by
\begin{equation}
\left(
\begin{array}{c}
V\\ U\\ \Theta \end{array}\right)(t,\x)=\left(
\begin{array}{c}
V^B+V^{CD}+V^R\\ U^B+U^{CD}+U^R\\ \Theta^B+\Theta^{CD}+\Theta^R
\end{array}\right)(t,\x)-\left(
\begin{array}{c}
v_*+v^*\\ u_*+u^*\\ \t_*+\t^* \end{array}\right), \label{(2.31)}
\end{equation}
where $(V^B, U^B, \Theta^B )(t,\x)$ is the subsonic BL-solution
defined in Lemma 2.1 (Case II) with the right state $(v_+, u_+,
\t_+)$ replaced by $(v_*, u_*, \theta_* )$, $(V^{CD}, U^{CD},
\Theta^{CD} )(t,\x)$ is the viscous contact wave   defined in
(\ref{(2.19)}) with the states $(v_-, u_-, \t_-)$ and $(v_+, u_+,
\t_+)$ replaced by $(v_*, u_*, \theta_* )$ and $(v^*, u^*, \theta^*
)$, respectively, and $(V^{R}, U^{R}, \Theta^{R} )(t,\x)$ is the
smoothed 3-rarefaction wave defined in (\ref{(2.26)}) with the left
state $(v_-, u_-, \t_-)$ replaced by $(v^*, u^*, \theta^* )$.

\vskip 2mm

\noindent{\bf Theorem 2.2} If $(v_-,u_-,\t_-) \in {\rm
BL\texttt{-}CD\texttt{-}R_3} (v_+,u_+,\t_+)$. Let $(V,U,\T)(t,x)$ be
the superposition of the BL-solution, the viscous contact wave and
the rarefaction wave defined in (\ref{(2.31)}). There exists a small
constant $\d_0$ such that if the BL-solution amplitude $\d^B$, the
contact discontinuity amplitude $\d^{CD}$ and the initial values
satisfy
$$
\d^B+\d^{CD}+\|(v_0-V_{0},u_0-U_{0},\t_0-\T_{0})\|_{1}\leq \d_0,
$$
then the inflow problem (\ref{(1.5)}) or the half space problem
(\ref{(1.5+)}) admits a unique global solution $(v,u,\t)(t,\x)$
satisfying
$$
\begin{array}{l}
\di (v-V,u-U,\t-\T)(t,\x)\in C([0,+\i),H^1({\bf R}^+)),\\[2mm]
\di (v-V)_\x(t,\x)\in L^2(0,+\i;L^2({\bf R}^+)),\\[2mm]
\di ((u-U)_\x,(\t-\T)_\x)(t,\x)\in L^2(0,+\i;H^1({\bf R}^+)),
\end{array}
$$
and
$$
\lim_{t\rightarrow +\i}\sup_{\x\in {\bf
R}^+}|(v-V,u-U,\t-\T)(t,\x)|=0.
$$

\section{Stability analysis}
\setcounter{equation}{0}

In this section we will prove our main stability results Theorem 2.1
and Theorem 2.2. We will focus on the proof of Theorem 2.2, i.e.,
the stability of the superposition wave. The proof of Theorem 2.1 is
almost same as Theorem 2.2 and we will omit it for brevity.

 Besides
the intrinsic properties of the BL-solution, the viscous contact
wave and the rarefaction wave in the stability analysis, the
interaction between the wave patterns should be dealt with carefully
in the stability analysis. Here we will use the elementary energy
methods to prove Theorem 2.2 by the classical continuum procedure.

Firstly we will reformulate the system of the superposition wave
$(V,U,\T)(t,\x)$ defined in (\ref{(2.31)}).

\subsection{Reformulation of the problem}

Recall the definition of the superposition wave $(V,U,\T)(t,\x)$
defined in (\ref{(2.31)}). Then we have
\begin{equation} \left\{
\begin{array}{l}
\di V_ {t}-\s_-V_{\xi}-U _{\xi}=0,\\
\di U_ {t}-\s_-U_{\xi}+P _{\xi}= \mu\big{(}\frac{U_{ \xi}}{V
}\big{)}_\xi+Q_1, \qquad
\qquad\qquad\qquad~~\,   \x>0,~t>0, \\
\di \f{R}{\g-1}(\T _{t}-\s_-\Theta_{\xi})+PU _{\xi}=\k( \f{\T _{
\xi}}{V }) _\x+\m  \f{U _{ \xi}^2}{V } +Q_2,\\
(V,U,\T)(t,0)=(v_-+V^{CD}-v_*,u_-+U^{CD}-u_*,\t_-+\T^{CD}-\t_*)(t,0),
\end{array}
\right. \label{(3.1)}
\end{equation}
where $ P=p(V, \Theta)=\f{R\Theta}{V}$, and the error terms
$Q_i~(i=1,2)$ are given by
\begin{equation}
\begin{array}{ll}
 \di Q_1=&\di
(P-P^B-P^{CD}-P^{R})_\x-\mu\left[(\frac{U_{\xi}}{V})_\x-(\frac{U^B_{\xi}}{V^B})_\x-(\f{U^{CD}_{\x}}{V^{CD}})_\xi\right]+\bar{Q}_1,\\[3mm]
\di Q_2=&\di
(PU_\x-P^BU^B_{\x}-P^{CD}U^{CD}_{\x}-P^RU^R_{\x})-\k\left[(\f{\T_{\xi}}{V
})_\x-(\f{\T^B_{\xi}}{V^B})_\x-(\f{\T^{CD}_{\xi}}{V^{CD}})
_\x\right]\\
&\di  -\m \left[\f{U_{\xi}^2}{V}-
\f{(U^B_{\xi})^2}{V^B}-\f{(U^{CD}_{\x})^2}{V^{CD}}\right]+\bar{Q}_2.
\end{array}
\label{(3.2)}
\end{equation}
and $\bar{Q}_i~(i=1,2)$ are the error terms defined in
(\ref{(2.23)}) to the viscous contact wave.

Due to the different propagation speeds of the BL-solution, the
viscous contact wave and the rarefaction wave, we can get the
following estimates of the error terms $Q_i(i=1,2)$:
\begin{equation}
\begin{array}{lll}
Q_1 &=&\di
O(1)\bigg[|(U^B_{\x},V^B_{\x},\Theta^B_{\x},U^B_{\x\x})||(V-V^B,\Theta-\Theta^B,V^{CD}_{\x},U^{CD}_{\x},V^{R}_{\x},U^{R}_{\x})|\\[2mm]
&&\di\qquad~~
+|(U^{CD}_{\x},V^{CD}_{\x},\Theta^{CD}_{\x},U^{CD}_{\x\x})||(V-V^{CD},\Theta-\Theta^{CD},V^{R}_{\x},U^{R}_{\x})|\\[2mm]
&&\di\qquad~~
+|(V^{R}_{\x},\Theta^{R}_{\x})||(V-V^{R},\Theta-\Theta^{R})|+|(U^{R}_{\x\x},U^{R}_{\x}V^{R}_{\x})|\bigg]+|\bar{Q}_1|\\[2mm]
&=&\di
O(1)(\d^B+\d^{CD})e^{-c(|\x|+t)}+O(1)(|U^{R}_{\x\x}|,|(U^{R}_{\x},V^{R}_{\x})|^2)+|\bar{Q}_1|,
\end{array}
\label{(3.3)}
\end{equation}
for some positive constant $c$ independent of $\x$ and $t$.
Similarly,
\begin{equation}
Q_2=O(1)(\d^B+\d^{CD})e^{-c(\x+t)}+O(1)(|\Theta^{R}_{\x\x}|,|(\Theta^{R}_{\x},V^{R}_{\x},U^{R}_{\x})|^2)+|\bar{Q}_2|.
\label{(3.4)}
\end{equation}

Denote the perturbation by
$$
(\p, \psi, \zeta)(t,\x)=(v, u, \t)(t,\x)-(V, U, \T)(t,\x),
$$
then we have the initial boundary value problem of the perturbation
$(\p,\psi,\z)(t,\x)$:
\begin{equation}
\left\{
\begin{array}{ll}
\di \p_ t-\s _-\p _\x-\psi _\x=0, &\x>0,~t>0,\\
\di \psi _t-\s _-\psi _\x+(p-P) _\x=\m(\f{u _\x}{v}-\f{U _\x}{V}) _\x-Q _1, &\x>0,~t>0,\\
\di \f{R}{\g-1}(\z _t-\s _-\z _\x)+(pu _\x-PU _\x)=\k(\f{\t
_\x}{v}-\f{\T _\x}{V}) _\x+\m(\f{u _\x^2}{v}-\f{U _\x^2}{V})-Q _2,&
\x>0,~t>0,\\[2mm]
\di
(\p,\psi,\z)(t,\x=0)=(v_--V,u_--U,\t_--\T)(t,\x=0),\\[2mm]
\di (\p,\psi,\z)(t=0,\x)=(\p _0,\psi _0,\z _0)(\x)\ra (0,0,0),~~{\rm
as}~~\x\ra+\i.
\end{array}
\right. \label{(3.5)}
\end{equation}

Since the local existence of the solution of (\ref{(3.5)}) is
well-known, we just state it and omit its proof for brevity.

Denote that
\begin{equation}
N(t)=\sup_{\tau\in[0,t]}\|(\p,\psi,\z)(\tau,\cdot)\|_{1}^2,\label{(3.6)}
\end{equation}
and define the solution space by
\begin{equation}
X_{\underline{m},\overline{m}}(0,T)=\left\{(\p,\psi,\z)(t,\x)\left|
\begin{array}{l}
\di (\p,\psi,\z)(t,\x)\in
 C([0,T];H^1({\bf R}^+)),\\[1mm]
\di (\psi_\x,\z_\x)\in L^2(0,T; H^1({\bf R}^+)),\\[1mm]
\di \p_\x\in L^2(0,T;L^2({\bf
 R}^+)),~~N(T)\leq
\overline{m}\\\di ~\inf_{[0,T]\times {\bf
R}^+}\{(V+\p),(\T+\z)\}(t,\x)\geq \underline{m}.
\end{array}
\right.\right\} \label{(3.7)}
\end{equation}
for some positive constants $\underline{m},\overline{m}$.

\vskip 2mm

\noindent{\bf Proposition 3.1 }(Local existence) Let
$(\p_0,\psi_0,\z_0)\in H^1({\bf R}^+)$. If
$\|(\p_0,\psi_0,\z_0)\|_1\leq \overline{m}$ and $\inf_{[0,T]\times
{\bf R}^+}\{(V+\p),(\T+\z)\}(t,\x)\geq \underline{m}$, then there
exist $\d_1$ and $t_0=t_0(\underline{m},\overline{m})>0$ such that
if the wave amplitude satisfies $\d^B+\d^{CD} < \d_1,$ then the half
space problem (\ref{(3.5)}) admits a unique solution
$(\p,\psi,\z)(t,\x)\in X_{\f {\underline{m}}2,
2\overline{m}}(0,t_0)$.

\vskip 2mm

 To prove Theorem 2.1, it
is sufficient to prove the following a priori estimate.

\vskip 2mm

\noindent{\bf Proposition 3.2 } (A priori estimate) Suppose that the
half space problem (\ref{(3.5)}) has a solution
$(\p,\psi,\z)(t,\x)\in X_{\f {\underline{m}}2, \v_0}[0,T]$ for a
suitably small constant $\v_0>0$. There exists a positive constant
$\d_2$, such that if the wave amplitude satisfies $\d^B+\d^{CD} <
\d_2,$ then the solution $(\p,\psi,\z)(t,\x)$ satisfy that for
$\forall~ t\in [0,T]$,
\begin{equation}
N(t)+\int_0^t\|\p_\x(\tau,\cdot)\|^2+\|(\psi_\xi,\z_\xi)(\tau,\cdot)\|_{1}^2d\tau
\leq C(N(0)+\d_2+\v^{\f18}),\label{(3.8)}
\end{equation}
where the positive constant $C$ is independent of $t$.

\subsection{Boundary estimates}
In this section we will obtain the boundary estimates needed in the
analysis below. From the definition of the viscous contact wave
(\ref{(2.19)}) and its property (\ref{(2.21)}), we have the
following estimates, which are very important in the boundary
estimates,
\begin{equation}
(v_*-V^{CD},u_*-U^{CD},\t_*-\T^{CD})(t,\x)=O(1)\d^{CD}
e^{-\f{c_0(\x+\s_-t)^2}{1+t}},~~ {\rm as}~~ |\xi+\s_-t|\ra \i.
\label{(3.9)}
\end{equation}
So on the boundary $\xi=0$,
\begin{equation}
\begin{array}{ll}
\di (\p,\psi,\z)(t,\x=0)&\di =(v_--V,u_--U,\t_--\T)(t,\x=0)\\
&\di =(v_*-V^{CD},u_*-U^{CD},\t_*-\T^{CD})(t,\x=0)\\
&\di =O(1)\d^{CD} e^{-\f{c_0(\s_-t)^2}{1+t}}\\
&\di =O(1)\d^{CD} e^{-c_1t},~~ {\rm as}~~ t\ra +\i.
\end{array}
\label{(3.10)}
\end{equation}
Thus we have the following lemma.

\vskip 2mm

\noindent{\bf Lemma 3.1}(Boundary estimates) There exists the
positive constant $C$ such that for any $t>0$,
\begin{equation}
\begin{array}{l}
\di \int_0^t|(\p,\psi,\z)|^2(\tau,\x=0)d\tau\leq C(\d^{CD})^2,\\[5mm]
\di \int_0^t[\m(\f{u_\xi}{v}-\f{U_\xi}{V})\psi](\tau,\x=0)d\tau\leq
\nu\int_0^t(\|\psi_{\x\x}\|^2+\|\psi_\x\|^2)d\tau +C_\nu (\d^{CD})^2,\\[5mm]
\di
\int_0^t[\k(\f{\t_\xi}{v}-\f{\T_\xi}{V})\f{\z}{\t}](\tau,\x=0)d\tau\leq
\nu\int_0^t(\|\z_{\x\x}\|^2+\|\z_\x\|^2)d\tau +C_\nu (\d^{CD})^2,\\[5mm]
\di \int_0^t [-\f{\m
\s_-}{2}(\f{\tilde{v}_{\x}}{\tilde{v}})^2+\psi\f{\tilde{v}_\tau}{\tilde{v}}](\tau,\x=0)d\tau\leq
\nu\int_0^t\|\psi_{\x\x}\|^2d\tau+C_\nu \int_0^t
\|\psi_\x\|^2d\tau+C(\d^{CD})^2,
\end{array}
\label{(3.11)}
\end{equation}
where $\di \tilde{v}=\f{v}{V}$, $\nu$ is a positive small constant
to be determined later and $C_\nu$ is a positive constant depending
on $\nu$.

\vskip 2mm

\noindent{\bf Proof:} The proof of $(\ref{(3.11)})_1$ is a direct
consequence of (\ref{(3.10)}).

Now we prove $(\ref{(3.11)})_2$.
$$
\begin{array}{ll}
&\di \int_0^t[\m(\f{u_\xi}{v}-\f{U_\xi}{V})\psi](\tau,\x=0)d\tau\\[3mm]
 &\di \leq
C\int_0^t |\psi|(|\psi_\x|+|U_\x||\p|)(\tau,\x=0) d\tau\\[3mm]
&\di \leq C\int_0^t |\psi_\x||\psi|(\tau,\x=0)d\tau
+C\int_0^t(|\psi|^2+|\p|^2)(\tau,\x=0)d\tau\\[4mm]
&\di \leq C\int_0^t \sup_{\x\in[0,+\i)}|\psi_\x(\tau,\x)|\cdot|\psi(t,\x=0)| d\tau+C(\d^{CD})^2\\[4mm]
&\di \leq
C\int_0^t\|\psi_\x(\tau,\cdot)\|^{\f12}\cdot\|\psi_{\x\x}(\tau,\cdot)\|^{\f12}\cdot|\psi(\tau,\x=0)|d\tau+C(\d^{CD})^2\\[4mm]
&\di \leq \nu
\int_0^t(\|\psi_{\x\x}(\tau,\cdot)\|^2+\|\psi_\x(\tau,\cdot)\|^2)d\tau+C_\nu\int_0^t|\psi(\tau,\x=0)|^2d\tau+C(\d^{CD})^2\\[4mm]
&\di \leq \nu
\int_0^t(\|\psi_{\x\x}(\tau,\cdot)\|^2+\|\psi_\x(\tau,\cdot)\|^2)d\tau+C_\nu(\d^{CD})^2.
\end{array}
$$
So the proof of $(\ref{(3.11)})_2$ is completed. Similarly, we can
obtain $(\ref{(3.11)})_3$.

Then we will verify the inequality $(\ref{(3.11)})_4$. Notice that
\begin{equation}
\f{\tilde{v}_{\x}}{\tilde{v}}=\f{v_\x}{v}-\f{V_\x}{V}=\f{\p_\x}{v}-\f{V_\x\p}{vV},\label{(3.11+)}
\end{equation}
and
$$
\begin{array}{ll}
\di \f{\tilde{v}_{t}}{\tilde{v}}&\di =\f{v_t}{v}-\f{V_t}{V}=\f{\s_-v_\x+u_\x}{v}-\f{\s_-V_\x+U_\x}{V}\\[2mm]
&\di =\s_-\f{\tilde{v}_{\x}}{\tilde{v}}+(\f{u_\x}{v}-\f{U_\x}{V}).
\end{array}
$$
So we have
$$
\begin{array}{ll}
&\di \int_0^t [-\f{\m
\s_-}{2}(\f{\tilde{v}_{\x}}{\tilde{v}})^2+\psi\f{\tilde{v}_\tau}{\tilde{v}}](\tau,\x=0)d\tau\\[4mm]
&\di~ \leq C\int_0^t
(|\p_\x|^2+|\psi_\x|^2+|\p|^2+|\psi|^2)(\tau,\x=0)d\tau\\
 &\di ~\leq C\int_0^t
(|\p_\tau|^2+|\psi_\x|^2+|\p|^2+|\psi|^2)(\tau,\x=0)d\tau\\[4mm]
&\di ~\leq C\int_0^t\sup_{\x\in[0,+\i)}|\psi_\x(\tau,\xi)|^2d\tau+C(\d^{CD})^2\\[4mm]
&\di ~\leq
C\int_0^t\|\psi_\x(\tau,\cdot)\|\|\psi_{\x\x}(\tau,\cdot)\|d\tau+C(\d^{CD})^2\\[4mm]
&\di ~\leq \nu
\int_0^t\|\psi_{\x\x}(\tau,\cdot)\|^2d\tau+C_\nu\int_0^t\|\psi_\x(\tau,\cdot)\|^2d\tau+C(\d^{CD})^2,
\end{array}
$$
where in the second inequality we have used the fact
$$
\p_\x=\f{\p_\tau-\psi_\x}{\s_-}.
$$
 Now Lemma 3.1
is proved.

\subsection{Energy estimates}

In this section we will prove the a priori estimate in Proposition
3.2. Firstly we have the following Lemma:

\vskip 2mm

\noindent{\bf Lemma 3.2} There exist a constant $C>0$ such that if
the wave amplitudes $\d^B$, $\d^{CD}$ and the constants $\v$, $\v_0$
are small enough, then we have $\forall~ t\in [0,T]$,
\begin{equation}
\begin{array}{ll}
&\di
\|(\p,\psi,\z,\phi_\x)(t,\cdot)\|^2+\int_0^t\|(\p_\x,\psi_\x,\z_\x)(\tau,\cdot)\|^2d\tau+\int_0^t\int_{\mathbf{R}^+}U^R_\x(\p^2+\z^2)d\x d\tau\\[3mm]
\leq &\di  C\|(\p_0,\psi_0,\z_0,\p_{0\x})\|^2 \di
+C(\d^B+\d^{CD}+\v^{\f18})\bigg[\int_0^t(1+\tau)^{-\f{13}{12}}\|(\p,\psi,\z)(\cdot,\tau)\|^2d\tau+1\bigg]\\[4mm]
&\di
+C\nu\int_0^t\|(\psi_{\x\x},\z_{\x\x})(\tau,\cdot)\|^2d\tau+C\d^{CD}
\int_0^t \int_{{\bf R}^+}
(1+\tau)^{-1}e^{\f{-c_0(\x+\s_-\tau)^2}{1+\tau}}|(\p,\z)|^2d\x
d\tau.
\end{array}\label{(3.12)}
\end{equation}
{\bf Proof:}
Let
$$
\Phi(z)=z-1-\ln z.
$$
Similar in \cite{[Huang-Xin-Yang]}, we can get the following
estimate
\begin{equation}
\begin{array}{ll}
& \di
I_{1t}(t,\x)+H _{1\x}(t,\x)+\mu\f{\T\psi_ \x^2}{v\t}+\k\f{\T\z _\x^2}{v\t^2}+PU^R_{\x}\left[\Phi(\f{\t V}{v\T})+\g\Phi(\f{v}{V})\right]\\[3mm]
&\di =Q_3-Q_1\psi-Q_2\f{\z}{\t},
\end{array}
\label{(3.13)}
\end{equation}
where
\begin{equation}
I_1(t,\x)=R\T\Phi(\f{v}{V})+\f{\psi^2}{2}+\f{R\T}{\g-1}\Phi(\f{\t}{\T}),
\label{(3.14)}
\end{equation}
\begin{equation}
H_1(t,\x)=-\s_-I_1(t,\x)+(p-P)\psi-\m(\f{u _\x}{v}-\f{U_
\x}{V})\psi-\k(\f{\t _\x}{v}-\f{\T _\x}{V})\f{\z}{\t},
\label{(3.15)}
\end{equation}
\begin{equation}
\begin{array}{ll}
Q_3=&\di -P(U^B_{\x}+U^{CD}_{\x})\left[\Phi(\f{\t
V}{v\T})+\g\Phi(\f{v}{V})\right]+\left[\k(\f{\T_\x}{V})_\x+\mu\f{U_\x^2}{V}+Q_2\right]\Big[(\g-1)\Phi(\f{v}{V})\\[3mm]
&\di+\Phi(\f{\t}{\T})-\f{\z^2}{\t\T}\Big] - \m(\f{1}{v}-\f{1}{V})U_
\x\psi_\x +\m(\f1v-\f1V)U_\x^2\f{\z}{\t}+2\m\f{\z\psi_\x
U_\x}{v\t}+\k
\f{\T_\x\z_\x\z}{v\t^2}\\[3mm]
&\di
-\k(\f{1}{v}-\f{1}{V})\f{\T\T_\x\z_\x}{\t^2}+\k(\f{1}{v}-\f{1}{V})\f{\z\T_\x^2}{\t^2}
\end{array}
\label{(3.16)}
\end{equation}
Note that
$$
\Phi(1)=\Phi^\prime(1)=0,\qquad \Phi^{\prime\prime}(z)=z^{-2}>0.
$$
So there exists a positive constant $C$  such that
$$
C^{-1}(z-1)^2\leq \Phi(z)\leq C(z-1)^2,
$$
if $z$ is near 1.

Using the a priori assumptions $N(T)\leq \v_0$ for suitably small
constant $\v_0$, we can get
\begin{equation}
C^{-1}|\p|^2\leq \Phi(\f{v}{V})\leq C|\p|^2,\qquad C^{-1}|\z|^2\leq
\Phi(\f{\t}{\T})\leq C|\z|^2, \label{(3.17)}
\end{equation}
and
\begin{equation}
C^{-1}|(\p,\z)|^2\leq\Phi(\f{\t V}{v\T})+\g\Phi(\f{v}{V})\leq
C|(\p,\z)|^2. \label{(3.18)}
\end{equation}
Substituting (\ref{(3.17)}) and (\ref{(3.18)}) into (\ref{(3.16)})
 and using Cauchy inequality imply
\begin{equation}
\begin{array}{ll}
\di Q_3\leq&\di  \f{\mu\T\psi_ \x^2}{4v\t}+\f{\k\T\z
 _\x^2}{4v\t^2}+O(1)\Big[|(V^B_{\x},U^B_{\x},\T^B_{\x},\T^B_{\x\x})|+(|\T^{CD}_{\x}|^2,|\T^{CD}_{\x\x}|)\\[3mm]
 &\di +(|(V^R_{\x},U^R_{\x},\T^R_{\x})|^2,|\T^R_{\x\x}|)+|Q_2|\Big](\p^2+\z^2)
\end{array}
\label{(3.19)}
\end{equation}
By the fact
\begin{equation}
|f(\x)|=|f(0)+\int^\x_0 f_y dy|\leq |f(0)|+ \sqrt{\x }~ \|f_\x\|,
\label{(3.20)}
\end{equation}
we have
\begin{equation}
\begin{array}{ll}
&\di \int_0^t \int_{R_+}
|(V^B_{\x},U^B_{\x},\T^B_{\x},\T^B_{\x\x})|(\p^2+\z^2)d\x d\tau\\
&\di \leq C\d^B \int_0^t \int_{R_+} e^{-c\x}\left(
|(\phi,\z)|^2(\tau,0)+\x\|(\p_\x,\z_\x)\|^2\right)d\x d\tau\\
&\di \leq C\d^B \int_0^t |(\phi,\z)|^2(\tau,0) d\tau+C\d^B\int_0^t
\|(\p_\x,\z_\x)\|^2 d\tau\\
&\di \leq C\d^B(\d^{CD})^2+C\d^B\int_0^t \|(\p_\x,\z_\x)\|^2 d\tau.
\end{array}
\label{(3.21)}
\end{equation}
By the properties of the viscous contact wave, we can obtain
\begin{equation}
\begin{array}{ll}
\di \int_0^t\int_{R_+}
(|\T^{CD}_{\x}|^2,|\T^{CD}_{\x\x}|)(\p^2+\z^2)d\x&\di\leq
   C \d^{CD} \int_0^t\int_{{\bf
R}^+}(1+\tau)^{-1}e^{-\f{c_0(\x+\s_-\tau)^2}{1+\tau}}|(\p,\z)|^2d\x.
\end{array}
\label{(3.22)}
\end{equation}
 Using the definition of the
approximate rarefaction wave, we have
\begin{equation}
\begin{array}{ll}
& \di \int_0^t\int_{R_+}
(|(V^R_{\x},U^R_{\x},\T^R_{\x})|^2,|\T^R_{\x\x}|)(\p^2+\z^2)d\x\\
& \di\leq \int_0^t
(\|(V^R_{\x},U^R_{\x},\T^R_{\x})\|^2 +\|\T^R_{\x\x}\|_{L^1} )\|(\p,\z)\|^2_{L^\infty}d\tau\\
&\di \leq C\v^ {\f18}\int_0^t
(1+\tau)^ {-\f{13}{16}}\|(\p,\z)\|\|(\p_\x,\z_\x)\|d\tau\\
&\di \leq C\v^ {\f18}\int_0^t (1+\tau)^
{-\f{13}{8}}\|(\p,\z)\|^2d\tau+ C\v^ {\f18}\int_0^t \|(\p_\x,
\z_\x)\|^2d\tau
\end{array}\label{(3.23)}
\end{equation}
where in the second inequality we have used
$$
\|(\Theta^{R}_{\x},V^{R}_{\x},U^{R}_{\x})\|^2 \leq
C\varepsilon^{\f18}(1+t)^{-\f{7}{8}}, $$ and
$$
\|(\Theta^{R}_{\x\x},V^{R}_{\x\x},U^{R}_{\x\x})\|_{L^1_\x}\leq
C\v^{\f18}(1+t)^{-\f{13}{16}},
$$
if we let $q\geq 16$ in Lemma 2.3.

Now we estimate the terms $Q_1\psi$, $Q_2\f{\z}{\t}$ on the
right-hand side of (\ref{(3.13)}) and the term $|Q_2|(\p^2+\z^2)$ on
the right-hand side of (\ref{(3.19)}). Due to the estimation of
$Q_1$ in (\ref{(3.3)}), we have
\begin{equation}
\begin{array}{ll}
\di \int_0^t \int_{\mathbf{R}^+}|Q_{1}\psi|d\x d\tau  &\di \leq C\int_0^t \|\psi\|_{L^\i_\x}\|Q_1\|_{L^1_\x} d\tau\\
 &\di \leq C\int_0^t \|\psi\|^{\f12}\|\psi_\x\|^{\f12}\bigg[\d^B e^{-c\tau}+\d^{CD}(1+\tau)^{-1}+C\v^{\f18}(1+\tau)^{-\f{13}{16}}\bigg]d\tau\\
 &\di \leq C(\d^B+\d^{CD}+\v^{\f18})\bigg[\int_0^t\|\psi_\x\|^2 d\tau+\int_0^t\|\psi\|^{\f23}(1+\tau)^{-\f{13}{12}}
 d\tau\bigg]\\[3mm]
 &\di \leq C(\d^B+\d^{CD}+\v^{\f18})\bigg[\int_0^t\|\psi_\x\|^2
 d\tau+\int_{0}^t\|\psi\|^2 (1+\tau)^{-\f{13}{12}} d\tau+1\bigg].
\end{array}\label{(3.24)}
\end{equation}

Similarly we can calculate the term $Q_2\f{\z}{\t}$ and
$|Q_2|(\p^2+\z^2)$.

Integrating (\ref{(3.14)}) over ${\bf R}^+\times [0,t]$ and using
the boundary estimates in Lemma 3.1, we can obtain
\begin{equation}
\begin{array}{ll}
&\di
\|(\p,\psi,\z)(t,\cdot)\|^2+\int_0^t\|(\psi_\x,\z_\x)(\tau,\cdot)\|^2d\tau+\int_0^t\int_{\mathbf{R}^+}U^R_\x(\p^2+\z^2)d\x d\tau\\
\leq &\di
C\|(\p_0,\psi_0,\z_0)\|^2+C(\d^B+\d^{CD}+\v^{\f18})\bigg[\int_{0}^t(1+\tau)^{-\f{13}{12}}\|(\p,\psi,\z)\|^2  d\tau+1\bigg]\\
&\di
+C\nu\int_0^t\|(\psi_{\x\x},\z_{\x\x})(\tau,\cdot)\|^2d\tau+C\d^{CD}\int_0^t\int_{{\bf
R}^+}(1+\tau)^{-1}e^{-\f{c_0(\x+\s_-\tau)^2}{1+\tau}}|(\p,\z)|^2d\x
d\tau.
\end{array}\label{(3.25)}
\end{equation}

Now we estimate $\|\p_\x\|^2$. Let
$$
\tilde v=\f{v}{V}.
$$
From the system $(\ref{(3.5)})_2$, we have
$$
\m(\f{\tilde{v}_\x}{\tilde v})_t-\m \s_- (\f{\tilde{v}_\x}{\tilde
v})_\x-\psi_t+\s_-\psi_\x-(p-P)_\x-Q_1=0.
$$
Multiplying the above equation by $\f{\tilde{v}_\x}{\tilde v}$ and
noticing that
$$
-(p-P)_\x=\f{R\t}{v}\f{\tilde{v}_\x}{\tilde
v}-\f{R\z_\x}{v}+(p-P)\f{V_\x}{V}-R\T_\x(\f1V-\f1v),
$$
we can get
$$
\begin{array}{ll}
&\di \left[\f{\m}{2}(\f{\tilde{v}_\x}{\tilde
v})^2-\psi\f{\tilde{v}_\x}{\tilde v}\right]_t-\left[\f{\m
\s_-}{2}(\f{\tilde{v}_\x}{\tilde v})^2-\psi\f{\tilde{v}_t}{\tilde
v}\right]_\x+\f{R\t}{v}(\f{\tilde{v}_\x}{\tilde v})^2\\[5mm]
=&\di
\psi_\x(\f{u_\x}{v}-\f{U_\x}{V})+\left[\f{R\z_\x}{v}-(p-P)\f{V_\x}{V}+R\T_\x(\f1V-\f1v)-Q_1\right]\f{\tilde{v}_\x}{\tilde
v}.
\end{array}
$$

Integrating the above equality and using the boundary estimate
(\ref{(3.11)}), we obtain
\begin{equation}
\begin{array}{ll}
&\di \int_{{\bf R}^+}\left[\f{\m}{2}(\f{\tilde{v}_\x}{\tilde
v})^2-\psi\f{\tilde{v}_\x}{\tilde
v}\right](t,\x)d\x+\int_0^t\int_{{\bf
R}^+}\f{R\t}{2v}(\f{\tilde{v}_\x}{\tilde v})^2d\x d\tau\\[4mm]
\leq &\di \int_{{\bf R}^+}\left[\f{\m}{2}(\f{\tilde{v}_\x}{\tilde
v})^2-\psi\f{\tilde{v}_\x}{\tilde
v}\right](0,\x)d\x+C(\d^{CD})^2+\nu\int_0^t\|\psi_{\x\x}(\tau,\cdot)\|^2d\tau
\\[4mm]
&\di+C\int_0^t\bigg[\|(\psi_\x,\z_\x)\|^2+\|Q_1\|^2\bigg]d\tau+C\int_0^t\int_{{\bf
R}^+}|(V_\x,U_\x,\T_\x)|^2|(\p,\z)|^2d\x d\tau.
\end{array}\label{(3.27)}
\end{equation}
Using the equality (\ref{(3.11+)}), we can get
\begin{equation}
C^{-1}(|\p_\x|^2-|V_\x\p|^2)\leq (\f{\tilde{v}_\x}{\tilde v})^2\leq
C(|\p_\x|^2+|V_\x\p|^2).\label{(3.28)}
\end{equation}
By the estimation of $Q_1$ in (\ref{(3.3)}), we have
\begin{equation}
\int_0^t\|Q_1\|^2d\tau\leq C(\d^B+\d^{CD}+\v^{\f18}).\label{(3.29)}
\end{equation}
Similar to (\ref{(3.21)})-(\ref{(3.23)}), we can compute the last
term in the right hand side of (\ref{(3.27)}). Thus we can obtain
\begin{equation}
\begin{array}{ll}
&\di \quad \|\phi_\x(t,\cdot)\|^2+\int_0^t\|\p_\x\|^2 d\tau \leq
C\|(\p_0,\psi_0,\p_{0\x})\|^2+C\|(\p,\psi)(t,\cdot)\|^2\\
&\di
+C\nu\int_0^t\|\psi_{\x\x}(\tau,\cdot)\|^2d\tau+C\d^{CD}\int_0^t\int_{{\bf
R}^+}(1+\tau)^{-1}e^{-\f{c_0(\x+\s_-\tau)^2}{1+\tau}}|(\p,\z)|^2d\x
d\tau\\
&\di+C\int_0^t\|(\psi_\x,\z_\x)\|^2d\tau
+C(\d^B+\d^{CD}+\v^{\f18})\bigg[\int_{0}^t(1+\tau)^{-\f{13}{12}}\|(\p,\psi,\z)\|^2
d\tau+1\bigg].
\end{array}
 \label{(3.30)}
\end{equation}
Multiplying the inequality (\ref{(3.25)}) by a large constant
$C_1>0$, and adding it to (\ref{(3.30)}), we complete the proof of
Lemma 3.2.

Now we derive the higher order estimates. Multiplying the equation
$(\ref{(3.5)})_2$ by $-\psi_{\x\x}$, we can get
\begin{equation}
(\f{\psi_\x^2}{2})_t-[\psi_t\psi_\x-\f{\s_- \psi_\x^2}{2}]_\x+\m
\f{\psi_{\x\x}^2}{v}=\m\f{\psi_\x}{v^2}v_\x
\psi_{\x\x}+\bigg\{(p-P)_\x-\m [U_\x
(\f{1}{v}-\f{1}{V})]_\x+Q_1\bigg\}\psi_{\x\x}.\label{(3.31)}
\end{equation}

From the boundary estimate
$$
\begin{array}{ll}
&\di \int_0^t[\psi_\tau\psi_\x-\f{\s_-
\psi_\x^2}{2}](\tau,\x=0)d\tau\\[3mm]
\leq &\di C
\int_0^t(|\psi_\x(\tau,0)|^2+|\psi_\tau(\tau,0)|^2)d\tau\\[3mm]
\leq &\di C \int_0^t \|\psi_\x\|\|\psi_{\x\x}\|d\tau+C(\d^{CD})^2\\[3mm]
\leq &\di \nu \int_0^t\|\psi_{\x\x}\|^2 d\tau+C_\nu \int_0^t
\|\psi_\x\|^2 d\tau +C(\d^{CD})^2,
\end{array}
$$
we can get the following inequality by integrating (\ref{(3.31)})
over ${\bf R}^+\times(0,t)$
\begin{equation}
\begin{array}{ll}
&\di \|\psi_\x\|^2(t) +\int_0^t\|\psi_{\x\x}\|^2d\tau \leq
C\|\psi_{0\x}\|^2+C(\d^{CD})^2+C\v^{\f15}\int_0^t(1+\tau)^{-\f32}\|(\p,\z)\|^2
d\tau \\
&\di \quad +C\int_0^t \|(\p_\x,\psi_\x,\z_\x)\|^2d\tau+C
(\d^{CD})^2\int_0^t \int_{{\bf R}^+}
(1+\tau)^{-1}e^{-\f{-c_0(\x+\s_-\tau)^2}{1+\tau}} |(\p,\z)|^2d\x
d\tau,
\end{array}\label{(3.32)}
\end{equation}
where we have used the following estimation
$$
\begin{array}{ll}
\di \int_0^t\int_{\mathbf{R}^+} |\p_\x||\psi_\x||\psi_{\x\x}|d\x
d\tau
&\di \leq C\int_0^t\|\p_\x\|\|\psi_{\x\x}\|\|\psi_\x\|_{L^\i_\x} d\tau\\
&\di \leq C\int_0^t
\|\p_\x\|\|\psi_{\x\x}\|^{\f32}\|\psi_\x\|^{\f12}
d\tau\\
&\di \leq \nu \int_0^t \|\psi_{\x\x}\|^2d\tau+C_\nu
\sup_{t}\|\p_\x\|^4\int_0^t \|\psi_\x\|^2d\tau\\
&\di \leq \nu \int_0^t \|\psi_{\x\x}\|^2d\tau+C_\nu ~\v_0^4\int_0^t
\|\psi_\x\|^2d\tau.
\end{array}
$$
Multiplying $(\ref{(3.5)})_3$ by $-\z_{\x\x}$, almost similar to the
estimates for $\|\psi_\x\|^2(t)$, we can obtain
\begin{equation}
\begin{array}{ll}
&\di \|\z_\x\|^2(t) +\int_0^t\|\z_{\x\x}\|^2d\tau \leq
C\|\z_{0\x}\|^2+C(\d^{CD})^2+C\v^{\f15}\int_0^t(1+\tau)^{-\f32}\|(\p,\z)\|^2
d\tau \\
&\di \quad +C\int_0^t \|(\p_\x,\psi_\x,\z_\x)\|^2d\tau+C
(\d^{CD})^2\int_0^t \int_{{\bf R}^+}
(1+\tau)^{-1}e^{-\f{-c_0(\x+\s_-\tau)^2}{1+\tau}} |(\p,\z)|^2d\x
d\tau.
\end{array}\label{(3.33)}
\end{equation}
Combining Lemma 3.2 and the higher order estimations (\ref{(3.32)})
and (\ref{(3.33)}), we have the following Lemma:

\vskip 2mm

\noindent{\bf Lemma 3.3} If the wave amplitudes $\d^B$, $\d^{CD}$
and the constants $\v$, $\v_0$ are small enough, then we have
$\forall~ t\in [0,T]$,
\begin{equation}
\begin{array}{ll}
&\di \|(\p,\psi,\z)(t,\cdot)\|_1^2+\int_0^t
\|\p_\x\|^2+\|(\psi_\x,\z_\x)\|_1^2 d\tau+\int_0^t\int_{\mathbf{R}^+}U^R_\x(\p^2+\z^2)d\x d\tau\\
\di \leq &\di
C\|(\p_0,\psi_0,\z_0)\|_1^2+C(\d^B+\d^{CD}+\v^{\f18})\bigg[\int_{0}^t(1+\tau)^{-\f{13}{12}}\|(\p,\psi,\z)\|^2
d\tau+1\bigg]\\
&\di +C \d^{CD}\int_0^t \int_{{\bf R}^+}
(1+\tau)^{-1}e^{-\f{-c_0(\x+\s_-\tau)^2}{1+\tau}} |(\p,\z)|^2d\x
d\tau.
\end{array}
\label{(3.34)}
\end{equation}

In order to close the estimate, we only need to control the last
term in (\ref{(3.34)}), which comes from the viscous contact wave.
So we will use the estimation on the heat kernel in
\cite{[Huang-Li-Matsumura]} and \cite{[Huang-Matsumura-Xin-2]}.

\vskip 2mm

 \noindent{\bf Lemma 3.4}  Suppose that $h(t,\x)$ satisfies
$$
h\in L^\i(0,T; L^2(\mathbf{R}^+)),~~h_\x\in L^2(0,T;
L^2(\mathbf{R}^+)),~~h_t-\s_-h_\x\in L^2(0,T; H^{-1}(\mathbf{R}^+)),
$$
Then
\begin{equation}
\begin{array}{ll}
&\di \int_0^t
\int_{\mathbf{R}^+}(1+\tau)^{-1}h^2e^{-\f{\b(\x+\s_-\tau)^2}{1+\tau}}d\x d\tau\\
&\di \leq C_\b\bigg[
\|h(0,\x)\|^2+\int_0^th^2(\tau,\x=0)d\tau+\int_0^t \|h_\x\|^2
d\tau+\int_0^t\langle h_t-\s_-h_\x,hg_\b^2\rangle_{H^1\times
H^{-1}}d\tau\bigg]
\end{array}\label{(3.35)}
\end{equation}
where
$$
g_\b(t,\x)=-(1+t)^{-\f12}\int_{\x+\s_-t}^{+\i} e^{-\f{\b
\eta^2}{1+t}}d\eta,
$$
and $\b>0$ is the constant to be determined.

The proof of Lemma 3.4 can be done similarly in
\cite{[Huang-Li-Matsumura]}. The only difference is that the space
we considered here is on the half line and the boundary terms should
be treated.

\vskip 2mm

 \noindent{\bf Lemma 3.5}  There exist a constant $C>0$ such that
if $\d^{CD}$ and $\v_0$ are small enough, then we have
\begin{equation}
\begin{array}{ll}
&\di \int_0^t\int_{{\bf R}^+}
\f{e^{-\f{c_0(\x+\s_-\tau)^2}{1+\tau}}}{1+\tau} |(\p,\psi,\z)|^2 d\x
d\tau\\[3mm]
&\di  \leq
C\bigg[(\d^B+\d^{CD}+\v^{\f18})+\|(\p_0,\psi_0,\z_0)\|^2+\|(\p,\psi,\z)(t,\cdot)\|^2\bigg]\\[2mm]
&~~\di +C\nu\int_0^t\|(\psi_{\x\x},\z_{\x\x})\|^2
d\tau+C\int_0^t\|(\p_\x,\psi_\x,\z_\x)\|^2d\tau+C\int_0^t
(1+\tau)^{-\f{13}{12}}\|(\p,\psi)\|^2 d\tau.
\end{array}\label{(3.36)}
\end{equation}
{\bf Proof:}
  From the equation $(\ref{(3.5)})_2$ and the fact
$p-P=\f{R\z-P\p}{v}$, we have
$$
\psi_t-s_-\psi_\x+(\f{R\z-P\p}{v})_\x=\m(\f{u_\x}{v}-\f{U_\x}{V})_\x-Q_1.
$$
Then we get
\begin{equation}
(R\z-P\p)_\x=\f{R\z-P\p}{v}(V_\x+\p_\x)-v(\psi_t-\s_-\psi_\x)+\m
v(\f{u_\x}{v}-\f{U_\x}{V})_\x-vQ_1.\label{(3.37)}
\end{equation}

Let
$$
G_\a(t,\x)=-(1+t)^{-1}\int_{\x+\s_-t}^{+\i} e^{-\f{\a
\eta^2}{1+t}}d\eta,
$$
where $\a$ is a positive constant to be determined later.
Multiplying the equation (\ref{(3.37)}) by $G_\a(R\z-P\p)$ gives
\begin{equation}
\begin{array}{ll}
&\di
\left[\f{G_\a(R\z-P\p)^2}{2}\right]_\x-(G_\a)_\x\f{(R\z-P\p)^2}{2}\\[2mm]
=&\di \f{G_\a(R\z-P\p)^2}{v}(V_\x+\p_\x)-G_\a v(R\z-P\p)(\psi_t-\s_-\psi_\x)\\[2mm]
&\di +\m G_\a v(R\z-P\p)(\f{u_\x}{v}-\f{U_\x}{V})_\x-G_\a
v(R\z-P\p)Q_1.
\end{array}\label{(3.38)}
\end{equation}
Note that
\begin{equation}
\begin{array}{ll}
\di -G_\a v(R\z-P\p)(\psi_t-\s_-\psi_\x)=-[G_\a
v(R\z-P\p)\psi]_t+[G_\a
v(R\z-P\p)\psi]_\x\\[3mm]
\di\quad +[(G_\a v)_t-\s_-(G_\a v)_\x](R\z-P\p)\psi+G_\a
v\psi[(R\z-P\p)_t-\s_-(R\z-P\p)_\x],
\end{array}\label{(3.39)}
\end{equation}
and
\begin{equation}
\begin{array}{ll}
&\di
(R\z-p_+\p)_t-\s_-(R\z-P\p)_\x\\[2mm]
&\di=(R\z_t-Rs_-\z_\x)-(P_t-\s_-P_\x)\p-P(\p_t-\s_-\p_\x)\\[2mm]
&\di =-\g
P\psi_\x+(\g-1)\bigg[-(p-P)(U_\x+\psi_\x)+\m(\f{u_\x^2}{v}-\f{U_\x^2}{V})+\k(\f{\t_\x}{v}-\f{\T_\x}{V})_\x
-Q_2\bigg]\\
&\di\quad -(P_t-\s_-P_\x)\p.
\end{array}\label{(3.40)}
\end{equation}
And using the equality
\begin{equation}
\di -G_\a v\g P \psi_\xi \psi =-[\g G_\a vP\f{\psi^2}{2}]_\x+\g v
P(G_\a)_\x\f{\psi^2}{2}+\g(vP)_\x\f{\psi^2}{2}, \label{(3.41)}
\end{equation}
 we can get
\begin{equation}
\f{e^{-\f{\a(\x+\s_-t)^2}{1+t}}}{2(1+t)}[(R\z-P\p)^2+\g
Pv\psi^2]=[G_\a v(R\z-P\p)\psi]_t+H_{2\x}(t,\x)+Q_4, \label{(3.42)}
\end{equation}
where
\begin{equation}
\begin{array}{ll}
\di H_2(t,\x)=&\di \f{G_\a(R\z-P\p)^2}{2}+\g G_\a vP\f{\psi^2}{2}-\s_- G_\a v(R\z-P\p)\psi\\[2mm]
&\di  -\k(\g-1)G_\a v \psi (\f{\t_\x}{v}-\f{\T_\x}{V})-\m G_\a v
(R\z-P \p)(\f{u_\x}{v}-\f{U_\x}{V}),
\end{array}
 \label{(3.43)}
\end{equation}
and
\begin{equation}
\begin{array}{ll}
\di Q_4=&\di -[(G_\a)_t-\s_-(G_\a)_\x]v(R\z-P\p)\psi -G_\a
u_\x(R\z-P\p)\psi \\[2mm]
&\di +(\g-1)G_\a
v\psi\left[(p-P)(U_\x+\psi_\x)-\m(\f{u_\x^2}{v}-\f{U_\x^2}{V})+Q_2\right]\\[2mm]
&\di+\m[G_\a v(R\z-P\p)]_\x(\f{u_\x}{v}-\f{U_\x}{V})+(\g-1)\k(G_\a v \psi)_\x(\f{\t_\x}{v}-\f{\T_\x}{V})\\[2mm]
&\di +G_\a v(R\z-P\p)Q_1+G_\a v\psi(P_t-\s_-P_\x)\p,
\end{array}
 \label{(3.44)}
\end{equation}
From the boundary estimate (\ref{(3.11)}), we have
\begin{equation}
\int_0^t H_2(\tau,\x=0) d\tau\leq C(\d^{CD})^2+\nu
\int_0^t\|(\psi_\x,\z_\x,\psi_{\x\x},\z_{\x\x})(\tau,\cdot)\|^2d\tau.
 \label{(3.45)}
\end{equation}
Note that
$$
\|G_\a(t,\cdot)\|_{L^\i}\leq C_\a (1+t)^{-\f12},
$$
thus integrating (\ref{(3.42)}) over $\mathbf{R}^+\times (0,t)$
gives
\begin{equation}
\begin{array}{ll}
&\di \int_0^t\int_{{\bf
R}^+}\f{e^{-\f{\a(\x+\s_-\tau)^2}{1+\tau}}}{1+\tau}[(R\z-P\p)^2+\psi^2]d\x d\tau\\[3mm]
&\di \leq
C\bigg[(\d^B+\d^{CD}+\v^{\f18})+\|\p_0,\psi_0,\z_0\|^2\bigg]+C(1+t)^{-1}\|(\p,\psi,\z)(t,\cdot)\|^2\\[2mm]
&\di +C\int_0^t
(1+\tau)^{-\f{13}{12}}\|(\p,\psi,\z)(\tau,\cdot)\|^2d\tau+C\int_0^t
\|(\p_\x,\psi_\x,\z_\x)(\tau,\cdot)\|^2 d\tau\\[2mm]
&\di +C\nu \int_0^t \|(\psi_{\x\x},\z_{\x\x})(\tau,\cdot)\|^2
d\tau+C\d^{CD}\int_0^t\int_{{\bf
R}^+}\f{e^{-\f{\a(\x+\s_-\tau)^2}{1+\tau}}}{1+\tau}|(\p,\z)|^2d\x
d\tau.
\end{array}\label{(3.46)}
\end{equation}

In order to get the desired estimate in Lemma 3.5, we must derive
the other similar estimates from the energy equation
$(\ref{(3.5)})_3$. Set
$$
h=R\z+(\g-1)p_+\p
$$
in Lemma 3.4. Thus we only need to compute the last term in
(\ref{(3.35)}). From the energy equation $(\ref{(3.5)})_3$, we have
\begin{equation}
h_t-\s_-h_\x=(P_t-\s_-P_\x)\p-(p-P)u_\x+\k(\f{\t_\x}{v}-\f{U_\x}{V})_\x+\m(\f{u_\x^2}{v}-\f{U_\x^2}{V})-Q_2,\label{(3.47)}
\end{equation}
Thus
\begin{equation}
\begin{array}{ll}
&\di \int_0^t\langle h_t-\s_-h_\x,hg_\b^2\rangle_{H^1\times
H^{-1}}d\tau\\
=&\di\int_0^t\int_{\mathbf{R}^+}[(P_t-\s_-P_\x)\p-(p-P)U_\x]hg_\b^2d\x
d\tau+\int_0^t\int_{\mathbf{R}^+}(p-P)\psi_\x hg_\b^2d\x d\tau\\
&\di+\int_0^t[\k(\f{\t_\x}{v}-\f{\T_\x}{V}) hg_\b^2](\tau,\x=0)
d\tau+
\int_0^t\int_{\mathbf{R}^+}\k(\f{\t_\x}{v}-\f{\T_\x}{V})(hg_\b^2)_\x
d\x d\tau\\
&\di
+\int_0^t\int_{\mathbf{R}^+}\m(\f{u^2_\x}{v}-\f{U^2_\x}{V})hg_\b^2
d\x d\tau+\int_0^t\int_{\mathbf{R}^+}Q_2 hg_\b^2 d\x d\tau\\
:=&\di \sum_{i=1}^6 J_i.
\end{array}\label{(3.48)}
\end{equation}
Note that
$$
\|g_\b(t,\cdot)\|_{L^\i}\leq C_\b,
$$
we can estimate $J_i(i=1,3,4,5,6)$ directly. In order to estimate
$J_2$, from the mass equation $(\ref{(3.5)})_1$, we have
$$
\begin{array}{ll}
&\di
(p-P)\psi_\x hg_\b^2\\[2mm]
=&\di \f{(\g-1)h-\g P\p}{v}hg_\b^2(\p_t-\s_-\p_\x)\\[2mm]
=&\di\f{(\g-1)h^2 g_\b^2}{v}(\p_t-\s_-\p_\x)-\f{\g Phg_\b^2}{2v}[(\p^2)_t-\s_-(\p^2)_\x]\\[2mm]
=&\di\Bigg[ \f{2(\g-1)\p h^2g_\b^2-\g Ph
\p^2g_\b^2}{2v}\Bigg]_t-\s_-\Bigg[ \f{2(\g-1)\p h^2g_\b^2-\g Ph
\p^2g_\b^2}{2v}\Bigg]_\x
\\[3mm]
&\di +\f{\g
Ph\p^2-2(\g-1)h^2\p}{v}g_\b[(g_\b)_t-\s_-(g_\b)_\x]-\f{\g
Ph\p^2-2(\g-1)h^2\p}{v^2}g_\b^2(v_t-\s_-v_\x)\\[3mm]
&\di+ \Bigg[\f{2(\g-1)g_\b^2\p h}{v}+\f{\g
Pg_\b^2\p^2}{2v}\Bigg](h_t-\s_-h_\x)+\f{\g
g_\b^2\p^2h}{2v}(P_t-\s_-P_\x)
\end{array}
$$

Now each term can be estimated directly, the detailed proof can be
seen in \cite{[Huang-Li-Matsumura]}. Remark that here we need to
compute the boundary terms. Therefore, taking the constant
$\b=\f{c_0}{2}$, we can get from Lemma 3.4 that
\begin{equation}
\begin{array}{ll}
&\di \int_0^t\int_{{\bf R}^+}
\f{e^{-\f{c_0(\x+\s_-\tau)^2}{1+\tau}}}{1+\tau} h^2 d\x d\tau \leq
C\bigg[(\d^B+\d^{CD}+\v^{\f18})+\|(\p_0,\psi_0,\z_0)\|^2+\|(\p,\psi,\z)(t,\cdot)\|^2\bigg]\\
&\quad \di +C\nu\int_0^t\|(\psi_{\x\x},\z_{\x\x})\|^2
d\tau+C\int_0^t\|(\p_\x,\psi_\x,\z_\x)\|^2d\tau+C\int_0^t
(1+\tau)^{-\f{13}{12}}\|(\p,\psi)\|^2 d\tau\\
&\quad \di +C(\d^{CD}+\v_0)\int_0^t \int_{{\bf
R}^+}(1+\tau)^{-1}e^{-\f{c_0(\x+\s_-\tau)^2}{1+\tau}}|(\p,\z)|^2d\x
d\tau.
\end{array}\label{(3.49)}
\end{equation}
Taking $\a=c_0$ in (\ref{(3.46)}) and combining the estimates
(\ref{(3.46)}) with (\ref{(3.49)}) yield the desired estimation in
Lemma 3.5 if we choose suitably small constants $\d^{CD}$ and
$\v_0$.

Now from Lemma 3.3 and Lemma 3.5, if the wave amplitude $\d^{CD}$
and the constant $\nu$ are suitably small, we can get
$$
\begin{array}{ll}
\di \|(\p,\psi,\z)(t,\cdot)\|_1^2  &\di +\int_0^t
\|\p_\x\|^2+\|(\psi_\x,\z_\x)\|_1^2 d\tau \leq
C\|(\p_0,\psi_0,\z_0)\|_1^2\\
&\di
+C(\d^B+\d^{CD}+\v^{\f18})\Bigg[\int_{0}^t(1+\tau)^{-\f{13}{12}}\|(\p,\psi,\z)\|^2
d\tau+1\Bigg].
\end{array}
$$
Finally, Gronwall inequality gives the a priori estimate in
Proposition 3.2:
$$
\|(\p,\psi,\z)(t,\cdot)\|_1^2+\int_0^t
\|\p_\x\|^2+\|(\psi_\x,\z_\x)\|_1^2 d\tau \leq
C\|(\p_0,\psi_0,\z_0)\|_1^2+C(\d^B+\d^{CD}+\v^{\f18}).
$$
Thus we complete the proof of Theorem 2.2.

The proof of Theorem 2.1 can be done along the same line as Theorem
2.2, we omit it for brevity.

\section*{Appendix: Proof of Lemma 2.1}
\renewcommand{\theequation}{A.\arabic{equation}}
\setcounter{equation}{0}

Now we give the rigorous proof of Lemma 2.1. Firstly from
(\ref{(2.9)}), $u_+>0$. Thus if $u_+\leq 0$, then there is no
solution to (\ref{(2.10)}) or (\ref{(2.13)}).  Now assume that
$u_+>0$. Then we can compute that the determinant of the matrix $J$
defined in (\ref{(2.12)})
\begin{equation}
\det J    =\frac{ R( u_+^2-R\gamma\theta_+) }{\k\m
(\g-1)}=\f{R^2\g\t_+(M_+^2-1)}{\k\mu(\g-1)}.\label{(A.1)}
\end{equation}

So we can divide it into three cases according to the sign of the
quantity $M_+^2 -1$.

\textbf{Case I (Supersonic):} ~$M_+>1$, then $\det J>0$. We can
easily know that $J$ has two positive eigenvalues. Thus the ODE
system (\ref{(2.13)}) has no solution.

\vskip 2mm

 \textbf{Case II (Transonic):} ~$M_+=1$,
then $\det J=0$. One of the eigenvalues of the matrix $J$ is zero,
the other one is positive. This case is a little subtle. Firstly we
can choose a nonsingular matrix $P$ such that $P^{-1}JP $ changes
into a standard form. For example, let
$$
P=\left(
\begin{array}{cc}
\f{\k(\g-1)^2}{R\k\mu+\k(\g-1)^2}&\f{R\k\g(\g-1)}{[R\k\mu+\k(\g-1)^2]u_+}\\[3mm]
-\f{\mu u_+}{\k(\g-1)}&1
\end{array}
\right).
$$
then
$$
P^{-1}JP = \left(
\begin{array}{cc}
\l_J &0\\ 0 & 0
\end{array}
\right):=\Lambda_J,
$$
where $\l_J$ is the positive eigenvalue of $J$ given by
$$
\l_J=\Bigg{(}\f{\g-1}{\mu}+\f{R}{\k(\g-1)}\Bigg{)}u_+>0.
$$
Let
\begin{eqnarray}
W=\left ( \begin{array}{c} W_1\\ W_2
\end{array}
\right):=P^{-1}\left ( \begin{array}{c} \bar{U}^B\\ \bar{\Theta}^B
\end{array}
\right),\label{(A.2)}
\end{eqnarray}
we have
\begin{eqnarray}
 W_{\x}=\Lambda_J W+G(W),\label{(A.3)}
\end{eqnarray}
where
$$
G(W)=P^{-1}F(PW),
$$
and
$$
F(PW)=\left(
\begin{array}{c}
F_1(PW)\\ F_2(PW)
\end{array}
\right).
$$
We can rewrite (\ref{(2.13)}) as
 \begin{eqnarray}
\left\{
\begin{array}{ll}
\di W_{1 \x} =\l_J W_1 + G_1(W_1, W_2 ),\\[1mm]
\di W_{2 \x} =G_2(W_1, W_2 ),
\end{array}
\right.\label{(A.4)}
\end{eqnarray}

Obviously, there exists a suitably small neighborhood
$\Omega_{\overline{\delta}_0}(0, 0)$ such that $(G_1,G_2) (W )$ is
analytic. And in this neighborhood, if $|W_1 |\ll|W_2 |$, then
\begin{eqnarray}
G_2(W )  &=& -\frac{   R^2  \gamma^2 \k(\gamma-1)(\gamma+1)  }{2[ R
\g\m+\k(\gamma-1)^2]^2u_+
  }W_2^2   +o\big{(} W _2^2\big{)} .\label{(A.5)}
\end{eqnarray}
From the geometric theory of the automatous ordinary differential
systems, we know that the equilibrium state (0,0) is a saddle-node
point to the system (\ref{(A.4)}). And (0,0) is an attractor whose
trajectory, denoted by $\Gamma$, is unique and tangent to $W_2$-axis
at (0,0). From the uniqueness of the attractor trajectory $\Gamma$,
we know that only when $( W_1, W_2)(0)\in \Gamma $, there exists a
solution to (\ref{(A.4)}), otherwise, there does not exist solution
to (\ref{(2.13)}). When $( W_1, W_2)(0)\in \Gamma $, the solution
$(W_1,W_2)(\x)$ satisfy that $| W _1(\x)|\ll| W _2(\x)|$ if $\x$ is
large enough, thus we have
\begin{eqnarray}
 -\sigma_1W_2^2\leq W_{2 \x}  \leq -\sigma _2 W_2^2,\label{(A.6)}
\end{eqnarray}
where $0<\sigma_1 <\sigma_2$ are two constants.

So we can get
\begin{eqnarray}
|(W_{1 }, W_{2 } )  (\x)| &\leq&  C   \frac{ \delta^B }{ 1+ \delta^B
 \x}  ~~~\x\in \mathbb{R}_+,\label{(A.7)}
\end{eqnarray}
where $\delta^B =|(W_{1 }, W_{2 } )(0)|=O(1)|(u_+-u_-,\t_+-\t_-)|$
is small enough.

From (\ref{(A.2)}), we can get the BL-solution $(U^B,\T^B)(\x)$ in
the transonic case ($M_+=1$) satisfy that
\begin{eqnarray}
&&   \frac{   \m u_+ }{ \k (\g-1) }(U^B-u_+) -(\T^B-\t_+)\cr\cr&=&
\int^\infty_\x\left[  -\frac{    u_+ (\bar{U}^B)^2}{ \k (\g-1) }
   +   \frac{(2-\g)u_+}{2\g\k } (\bar{U}^B)^2 +\f{R  }{\k(\g-1)}
   \bar{U}^B
\bar{\Theta}^B -\frac{ 1}{2\k} (\bar{U}^B)^3 \right]
d\x,\label{(A.8)}
\end{eqnarray}

\vskip 2mm
 \textbf{Case III (Subsonic):} ~$M_+<1$, then $\det J<0$.
One can see that $J$ has one positive and one negative eigenvalues.
Similar to Case II, we can choose a nonsingular matrix $P$ such that
$P^{-1}JP $ is in a standard form. Now we give the detailed
procedure for choosing the matrix $P$. Firstly, let
$$
P_1=\left ( \begin{array}{cc} 1&0\\a_1 &1
\end{array}
\right),
$$
where the constant $a_1$ is to be determined, then
$$
P_1^{-1}=\left ( \begin{array}{cc} 1&0\\-a_1 &1
\end{array}
\right),
$$
Assume that
\begin{eqnarray}
P_1^{-1} J P_1 = \left ( \begin{array}{cc}  m_{11 }& m_{12}
\\ 0 & m_{22}
\end{array}
\right):=\mathfrak{M} ,\label{(A.9)}
\end{eqnarray}
where the constants $m_{11}, m_{12}, m_{22}$ will be fixed when
$a_1$ is determined. From $m_{21}=0$ in (\ref{(A.9)}), we get a
equation of $a_1$:
$$
 \frac{R}{  \m}a_1^2+ \Bigg{(}\frac{(M_+^2\gamma-1)u_+}{ M_+^2\g   \m }
 -\frac{   R u_+ }{ \k(\g-1)  }\Bigg{)}a_1-\frac{u_+^2 }{ M_+^2  \gamma    \k  } =0 ,
$$
i.e.,
$$
  \Bigg{(}\frac { a_1 }{  u_+ }\Bigg{)}^2+ \Bigg{(}\frac{ M_+^2\gamma-1  }
  {    M_+ ^2R\gamma } -\frac{\m }{ \k(\g-1)  }\Bigg{)} \frac {   a_1  }{  u_+ }  -\frac{    \m }{ M_+^2 R \gamma    \k  }=0.
$$
Then we can solve the above equation to obtain
\begin{eqnarray}
 a_1=   c_1 u_+   <0~~~~ {\rm or} ~~ ~  a_1= c_2u_+>0
 ,\label{(A.10)}
\end{eqnarray}
where $c_1 < \min \Big{\{}0, -\frac{ M_+^2\gamma-1  }
  { M_+^2   R \gamma }\Big{\}},  \, c_2 >\max \Big{\{}\frac{\m}{ \k(\g-1)} , \frac{\m}{\k(\g-1)}-\frac{M_+^2\gamma-1}
  { M_+^2   R \gamma  }\Big{\}}>0 $ are the
solutions of the following equation
\begin{eqnarray}
  y^2+ \Bigg{(}\frac{ M_+^2\gamma-1  }
  {    M_+ ^2R\gamma } -\frac{\m }{ \k(\g-1) }\Bigg{)} y  -\frac{ \m }{ M_+^2 R \gamma    \k  }=0.\label{(A.11)}
\end{eqnarray}
Without loss of generality, we choose $a_1=c_2u_+$, then we can
compute that the matrix $\mathfrak{M}$ in (\ref{(A.9)})
$$
\mathfrak{M}= \left ( \begin{array}{cc}
  \l_J^1& \frac{R}{ \m}\\[2mm]   0 &\l_J^2
\end{array} \right)
$$
where $\l_J^1= \Big{(}\frac{ M_+^2 \gamma-1}{M_+^2 R  \gamma }+
c_2\Big{)}
  u_+>0$
and $\l_J^2=\Big{(}\frac{    R }{ \k(\g-1)  }- c_2\Big{)} u_+<0$ are
the two eigenvalues of the matrix $J$.

Then we can choose a matrix
$$
P_2=\left ( \begin{array}{cc} 1&a_2\\0 &1
\end{array}
\right),
$$
such that
\begin{eqnarray}
P_2^{-1} \mathfrak{M} P_2 = \left ( \begin{array}{cc}  \l_J^1&0
\\ 0 &\l_J^2
\end{array}
\right):=\Lambda_J ,\label{(A.12)}
\end{eqnarray}

Then we can get
\begin{eqnarray}
 a_2 =-\frac{R}{\m (\l_J^1-\l_J^2) } .\label{(A.13)}
\end{eqnarray}

Now we set
$$
P=P_1P_2=\left ( \begin{array}{cc}
  1& \f{a_2}{u_+}
\\  c_2u_+   &1+a_2c_2
\end{array}
\right)
$$
Then
$$
P^{-1}JP=\Lambda_J={\rm diag}\{\l_J^1,\l_J^2\}.
$$
Let
\begin{eqnarray}
W=\left ( \begin{array}{c} W_1\\ W_2
\end{array}
\right):=P^{-1}\left ( \begin{array}{c} \bar{U}^B\\ \bar{\Theta}^B
\end{array}
\right),\label{(A.14)}
\end{eqnarray}
we have
\begin{eqnarray}
 W_{\x}=\Lambda_J W+G(W),\label{(A.15)}
\end{eqnarray}
where
$$
G(W)=P^{-1}F(PW),
$$
and
$$
F(PW)=\left(
\begin{array}{c}
F_1(PW)\\ F_2(PW)
\end{array}
\right).
$$

We can rewrite (\ref{(A.15)}) as
 \begin{eqnarray}
\left\{
\begin{array}{ll}
\di W_{1 \x} =\lambda_J^1 W_1 + G_1(W_1, W_2 ),\\[1mm]
\di W_{2 \x} =\l_J^2 W_2+G_2(W_1, W_2 ),
\end{array}
\right.\label{(A.16)}
\end{eqnarray}

From above, one can easily know that $G_{1 }, G_{2}$ are analytic
with respect to $(W_1, W_2)$ near $(0, 0)$,
 then the equilibrium point $(0, 0)$ is the saddle point of (\ref{(A.16)}), i.e., in a suitably small neighborhood
$\Omega_{\overline{\delta}_0}(0, 0)$,  there exist    two opposite
attractor trajectories $\Gamma_{1 }, \Gamma_2$ tangent to $W_2$-axis
at $(0, 0)$. Let $\mathcal{M} =\Gamma_{1 }\cup \Gamma_2$, then
$\mathcal{M} $ is a  center-stable manifold. Only when $( W_1,
W_2)(0)\in \mathcal{M} $, there exists a solution the ODE system
(\ref{(A.16)}). In such case, there exist two positive constants
$\s_3, \s_4 $ which is close to $-\l_A^2$ such that
\begin{eqnarray}
  -\s_3 W_2\leq W_{2  }^\prime \leq -\s_4  W_2,  ~~~\x\in \mathbb{R}_+.\label{(A.17)}
\end{eqnarray}
 So we have that there exist positive constants $ c$ and $C$ such that
\begin{eqnarray}
  |( W_1, W_2)(\x)|\leq  C\d^B e^{-c \x}, ~~~\x\in
  \mathbb{R}_+,\label{(A.18)}
\end{eqnarray}
where $\d^B=|(W_1,W_2)(0)|=O(1)|(u_+-u_-,\t_+-\t_-)|$ is the
amplitude of the BL-solution.

The BL-solution $(U^B,\Theta^B)$ satisfies
\begin{eqnarray}
  &&( 1+a_2c_2u_+)(U^B-u_+)
-a_2(\T^B-\t_+) \cr\cr&= &\int^\infty_\x e^{-\l_A^2 \x}\bigg{\{}
\frac{1+a_2c_2}{\m}     (\bar{U}^B)^2+\f{R  }{\k(\g-1)} \bar{U}^B
\bar{\Theta}^B
 \cr&  &\qquad\qquad~~~-\frac{a_2    }{u_+ }\bigg{[}\bigg{(}
 \frac{Ru_+}{M_+^2\k \g }-\frac{u_+}{2\k}\bigg{)} (\bar{U}^B)^2 -\frac{ 1}{2\k} (\bar{U}^B)^3\bigg{]} \bigg{\}}
d\x.\label{(A.19)}
\end{eqnarray}
Now we complete the proof of Lemma 2.1.

\end{document}